\documentclass{amsart}

\usepackage{amssymb,dsfont,fullpage,units}

\newtheorem{lemma}{Lemma}
\newtheorem{theorem}{Theorem}

\newtheorem{proposition}{Proposition}
\theoremstyle{definition}
\newtheorem{condition}{Condition}
\newtheorem{example}{Example}
\theoremstyle{remark}
\newtheorem{remark}{Remark}

\def\R{\mathbb{R}}
\newcommand\dotp[1]{\langle #1 \rangle}

\def\wt{\widetilde}
\def\wh{\widehat}
\def\t{{\scriptscriptstyle\top}}
\def\F{{\operatorname{F}}}

\def\E{\mathbb{E}}
\def\Eh{\wh{\E}}

\def\Sig{\varSigma}
\def\Sigh{\wh{\Sig}}
\def\Sigl{\Sig_{\lambda}}
\def\Sighl{\Sigh_{\lambda}}

\def\xt{\tilde{x}}
\def\Sigt{\wt{\Sig}}
\def\Sigw{\Sig_w}
\def\Sighw{\Sigh_w}
\def\Sighlw{\Sigh_{\lambda,w}}

\def\dl{d_{1,\lambda}}
\def\dlt{\tilde{d}_{1,\lambda}}
\def\dll{d_{2,\lambda}}

\def\bh{\hat{\beta}}
\def\bb{\bar{\beta}}
\def\blh{\bh_{\lambda}}
\def\bl{\beta_{\lambda}}
\def\blb{\bar{\beta}_{\lambda}}
\def\bols{\bh_0}

\def\Deltal{\Delta_{\lambda}}

\def\bias{\operatorname{approx}}
\def\biasl{\operatorname{approx}_{\lambda}}
\def\noise{\operatorname{noise}}

\def\rhol{\rho_{\lambda}}
\def\snoise{\sigma}

\def\bbias{b}
\def\bbiasl{b_{\lambda}}

\def\epsI{\varepsilon_{\operatorname{rg}}}
\def\epsII{\varepsilon_{\operatorname{bs}}}
\def\epsIII{\varepsilon_{\operatorname{vr}}}
\def\deltaI{\delta_{\operatorname{s}}}
\def\deltaII{\delta_{\operatorname{f}}}

\DeclareMathOperator{\diag}{diag}
\DeclareMathOperator{\tr}{tr}
\DeclareMathOperator{\var}{var}

\newcommand{\ignore}[1]{}

\begin{document}

\title{Random design analysis of ridge regression}

\author[Hsu]{Daniel Hsu}
\address[D.~Hsu]{
  Department of Computer Science \\
  Columbia University \\
  450 Computer Science Building \\
  1214 Amsterdam Avenue, Mailcode: 0401 \\
  New York, NY 10027-7003}
\email[D.~Hsu]{djhsu@cs.columbia.edu}

\author[Kakade]{Sham M. Kakade}
\address[S.~M.~Kakade]{
  Microsoft Research \\
  One Memorial Drive \\
  Cambridge, MA, 02142}
\email[S.M.~Kakade]{skakade@microsoft.com}

\author[Zhang]{Tong Zhang}
\address[T.~Zhang]{
  Department of Statistics \\
  Rutgers University \\
  501 Hill Center \\
  110 Frelinghuysen Road \\
  Piscataway, NJ 08854}
  \email[T.~Zhang]{tzhang@stat.rutgers.edu}



\begin{abstract}
This work gives a simultaneous analysis of both the ordinary least squares
estimator and the ridge regression estimator in the random design setting
under mild assumptions on the covariate/response distributions.
In particular, the analysis provides sharp results on the ``out-of-sample''
prediction error, as opposed to the ``in-sample'' (fixed design) error.
The analysis also reveals the effect of errors in the estimated covariance
structure, as well as the effect of modeling errors, neither of which
effects are present in the fixed design setting.
The proofs of the main results are based on a simple decomposition lemma
combined with concentration inequalities for random vectors and matrices.
\end{abstract}

\maketitle

\section{Introduction}

In the random design setting for linear regression, we are provided
with samples of covariates and responses, 
$(x_1,y_1),(x_2,y_2),\dotsc,(x_n,y_n)$, which are
 sampled
independently from a population, where the $x_i$ are random vectors and the
$y_i$ are random variables.
Typically, these pairs are hypothesized to have the linear relationship
\[ y_i = \dotp{\beta,x_i} + \epsilon_i \]
for some linear function $\beta$ (though this hypothesis need not be true).
Here, the $\epsilon_i$ are error terms,
typically assumed to be normally distributed as $\mathcal{N}(0,\snoise^2)$.
The goal of estimation in this setting is to find coefficients $\bh$ based
on these $(x_i,y_i)$ pairs such that the expected prediction error on a new
draw $(x,y)$ from the population, measured as $\E[(\dotp{\bh,x} - y)^2]$,
is as small as possible.
This goal can also be interpreted as estimating $\beta$ with accuracy
measured under a particular norm. 

The random design setting stands in contrast to the fixed design setting,
where the covariates $x_1,x_2,\dotsc,x_n$ are fixed (\emph{i.e.},
deterministic), and only the responses $y_1,y_2,\dotsc,y_n$ treated as
random.
Thus, the covariance structure of the design points is completely known and
need not be estimated, which simplifies the analysis of standard
estimators.
However, the fixed design setting does not directly address out-of-sample
prediction, which is of primary concern in many applications;
for instance, in prediction problems, the estimator $\bh$ is computed from
an initial sample from the population, and the end-goal is to use $\bh$ as
a predictor of $y$ given $x$ where $(x,y)$ is a new draw from the
population.
A fixed design analysis only assesses the accuracy of $\bh$ on data already
seen, while a random design analysis is concerned with the predictive
performance on unseen data.

This work gives a detailed analysis of both the ordinary least squares
and \emph{ridge} estimators~\cite{Hoerl62} in the random design setting
that quantifies the essential differences between random and fixed
design.
In particular, the analysis reveals, through a simple decomposition:
\begin{itemize}
  \item the effect of errors in the estimated covariance structure;
    
  \item the effect of errors in the estimated covariance structure, as well
    as the effect of approximating the true regression function by a linear
    function in the case the model is misspecified;

  \item the effect of errors due to noise in the response.
\end{itemize}
Neither of the first two effects is present in the fixed design analysis of
ridge regression, and the random design analysis shows that the effect of
errors in the estimated covariance structure is minimal---essentially a
second-order effect as soon as the sample size is large enough.
The analysis also isolates the effect of approximation error in the main
terms of the estimation error bound so that the bound reduces to one that
scales only with the noise variance when the approximation error vanishes.

Another important feature of the analysis that distinguishes it from that
of previous work is that it applies to the ridge estimator with an
arbitrary setting of $\lambda \geq 0$.
The estimation error is given in terms of the spectrum of the second moment
of $x$ and the particular choice of $\lambda$---the dimension of the
covariate space does not enter explicitly except when $\lambda=0$.
When $\lambda = 0$, we immediately obtain an analysis of ordinary least
squares; we are not aware of any other random design analysis of the ridge
estimator with this characteristic.
More generally, the convergence rate can be optimized by appropriately
setting $\lambda$ based on assumptions about the spectrum.

Finally, while our analysis is based on an operator-theoretical approach
similar to that of \cite{SmaZou07} and \cite{CD07}, it relies on
probabilistic tail inequalities in a modular way that gives explicit
dependencies without additional boundedness assumptions other than those
assumed by the probabilistic bounds.

\medskip\noindent {\bf Outline.}
Section~\ref{section:preliminaries} discusses the model, preliminaries, and
related work.
Section~\ref{section:ridge} presents the main results on the excess mean
squared error of the ordinary least squares and ridge estimators under
random design and discusses the relationship to the standard fixed design
analysis.
Section~\ref{section:applications} discusses an application to accelerating
least squares computations on large data sets.
The proofs of the main results are given in
Section~\ref{section:ridge-proof}.

\section{Preliminaries}
\label{section:preliminaries}

\subsection{Notation}

Unless otherwise specified, all vectors in this work are assumed to live in
a finite dimensional inner product space with inner product
$\dotp{\cdot,\cdot}$.
The restriction to finite-dimensions is due to the probabilistic bounds
used in the proofs; the main results of this work can be extended to
(possibly infinite-dimensional) separable Hilbert spaces under mild
assumptions by using suitable infinite-dimensional generalizations of these
probabilistic bounds.
We denote the dimensionality of this space by $d$, but stress that our
results will not explicitly depend on $d$ except when considering the
special case of $\lambda=0$.
Let $\|\cdot\|_M$ for a self-adjoint positive definite linear operator
$M \succ 0$ denote the vector norm given by $\|v\|_M :=
\sqrt{\dotp{v,Mv}}$.
When $M$ is omitted, it is assumed to be the identity $I$, so $\|v\| =
\sqrt{\dotp{v,v}}$.
Let $u \otimes u$ denote the outer product of a vector $u$, which acts as
the rank-one linear operator $v \mapsto (u \otimes u)v = \dotp{v,u} u$.
For a linear operator $M$, let $\|M\|$ denote its spectral (operator) norm,
\emph{i.e.}, $\|M\| = \sup_{v \neq 0} \|Mv\| / \|v\|$, and let $\|M\|_\F$
denote its Frobenius norm, \emph{i.e.}, $\|M\|_\F = \sqrt{\tr(M^* M)}$.
If $M$ is self-adjoint, $\|M\|_\F = \sqrt{\tr(M^2)}$.
Let $\lambda_{\max}[M]$ and $\lambda_{\min}[M]$, respectively, denote the
largest and smallest eigenvalue of a self-adjoint linear operator $M$.

\subsection{Linear regression}

Let $x$ be a random vector, and let $y$ be a random variable.
Throughout, it is assumed that $x$ and $y$ have finite second moments
($\E[\|x\|^2] < \infty$ and $\E[y^2] < \infty$).
Let $\{v_j\}$ be the eigenvectors of
\begin{equation} \label{eq:Sig}
\Sig := \E[x \otimes x] ,
\end{equation}
so that they form an orthonormal basis.
The corresponding eigenvalues are
\[ \lambda_j := \dotp{v_j,\Sig v_j} = \E[\dotp{v_j,x}^2] . \]
It is without loss of generality that we assume all eigenvalues $\lambda_j$
are strictly positive, since otherwise we may restrict attention of all
vectors to a subspace in which the assumption holds.
Let $\beta$ achieve the minimum \emph{mean squared
error} over all linear functions, \emph{i.e.},
\[ \E[(\dotp{\beta,x} - y)^2]
= \min_w \left\{ \E[(\dotp{w,x} - y)^2] \right\} , \]
so that
\begin{equation} \label{eq:beta}
\beta := \sum_j \beta_j v_j
\quad \text{where} \quad
\beta_j := \frac{\E[\dotp{v_j,x}y]}{\E[\dotp{v_j,x}^2]}
.
\end{equation}
\ignore{
It is easy to check that $\beta$ achieves the minimum \emph{mean squared
error} over all linear functions, \emph{i.e.},
\[ \E[(\dotp{\beta,x} - y)^2]
= \min_w \left\{ \E[(\dotp{w,x} - y)^2] \right\} . \]

Indeed, since $\|v\|_\Sig^2 = \dotp{v,\E[x \otimes x]v} =
\E[\dotp{v,x}^2]$, 
}
We also have that the \emph{excess} mean squared error of
$w$ over the minimum is
\[ \E[(\dotp{w,x}-y)^2] - \E[(\dotp{\beta,x}-y)^2] = \|w-\beta\|_\Sig^2
\]
(see Proposition~\ref{proposition:regret}).

\subsection{The ridge and ordinary least squares estimators}

Let $(x_1,y_1), (x_2,y_2), \dotsc, (x_n,y_n)$ be independent copies of
$(x,y)$, and let $\Eh$ denote the empirical expectation with respect to
these $n$ copies, \emph{i.e.},
\begin{equation} \label{eq:Sigh}
 \Eh[f] := \frac1n \sum_{i=1}^n f(x_i,y_i)  \quad \quad
\Sigh := \Eh[x \otimes x]
= \frac1n \sum_{i=1}^n x_i \otimes x_i .
\end{equation}

Let $\blh$ denote the \emph{ridge estimator} with parameter $\lambda \geq
0$, defined as the minimizer of the $\lambda$-regularized empirical mean
squared error, \emph{i.e.},
\begin{equation} \label{eq:blh}
\blh := \arg\min_w
\left\{ \Eh[(\dotp{w,x} - y)^2] + \lambda \|w\|^2 \right\} .
\end{equation}
The special case with $\lambda = 0$ is the \emph{ordinary least squares
estimator}, which minimizes the empirical mean squared error.
These estimators are uniquely defined if and only if $\Sigh + \lambda I
\succ 0$ (a sufficient condition is $\lambda > 0$), in which case
\[ \blh = (\Sigh + \lambda I)^{-1} \Eh[xy] . \]

\subsection{Data model}
\label{section:datamodel}

We now specify the conditions on the random pair $(x,y)$ under
which the analysis applies.

\subsubsection{Covariate model}

We first define the following effective dimensions of the covariate $x$
based on the second moment operator $\Sig$ and the regularization level
$\lambda$:
\begin{equation} \label{eq:dim}
d_{p,\lambda} := \sum_j \left( \frac{\lambda_j}{\lambda_j + \lambda}
\right)^p
, \quad p \in \{1,2\}
.
\end{equation}
It will become apparent in the analysis that these dimensions govern the
sample size needed to ensure that $\Sig$ is estimated with sufficient
accuracy.
For technical reasons, we also use the quantity
\begin{equation} \label{eq:dlt}
\dlt := \max\{ \dl, 1 \}
\end{equation}
merely to simplify certain probability tail inequalities in the main result
in the peculiar case that $\lambda \to \infty$ (upon which $\dl \to 0$).
We remark that $\dll$ appears naturally arises in the standard fixed design
analysis of ridge regression (see Proposition~\ref{proposition:fixed}), and
that $\dl$ was also used by~\cite{Zhang05} and \cite{CD07} in their
random design analyses of (kernel) ridge regression.
It is easy to see that $\dll \leq \dl$, and that $d_{p,\lambda}$ is at most
the dimension $d$ of the inner product space (with equality iff $\lambda =
0$).

Our main
condition requires that the squared length of $(\Sig + \lambda
I)^{-1/2} x$ is never more than a constant factor greater than its
expectation (hence the name \emph{bounded statistical leverage}).
The linear mapping $x \mapsto (\Sig + \lambda I)^{-1/2} x$ is sometimes
called \emph{whitening} when $\lambda = 0$.
The reason for considering $\lambda > 0$, in which case we call the mapping
\emph{$\lambda$-whitening}, is that the expectation $\E[\|(\Sig + \lambda
I)^{-1/2}x\|^2]$ may only be small for sufficiently large $\lambda$, as
\[
\E[\|(\Sig + \lambda I)^{-1/2} x\|^2]
= \tr((\Sig + \lambda I)^{-1/2} \Sig (\Sig + \lambda I)^{-1/2})
= \sum_j \frac{\lambda_j}{\lambda_j + \lambda}
= \dl
.
\]
\begin{condition}[Bounded statistical leverage at $\lambda$]
\label{cond:leverage}
There exists finite $\rhol \geq 1$ such that, almost surely,
\[
\frac{\|(\Sig + \lambda I)^{-1/2} x\|}
{\sqrt{\E[\|(\Sig + \lambda I)^{-1/2} x\|^2]}}
= \frac{\|(\Sig + \lambda I)^{-1/2} x\|}{\sqrt{\dl}}
\leq \rhol
.
\]
\end{condition}
The hard ``almost sure'' bound in Condition~\ref{cond:leverage} may be
relaxed to moment conditions simply by using different probability tail
inequalities in the analysis.
We do not consider this relaxation for sake of simplicity.
We also remark that it is possible to replace
Condition~\ref{cond:leverage} with a subgaussian condition (specifically, a
requirement that every projection of $(\Sig + \lambda I)^{-1/2} x$ be
subgaussian), which can lead to a sharper deviation bound in certain cases.

\begin{remark}[Ordinary least squares]
\label{remark:leverage-ols}
If $\lambda = 0$, then Condition~\ref{cond:leverage} reduces to the
requirement that there exists a finite $\rho_0\geq 1$ such that, almost
surely,
\[
\frac{\|\Sig^{-1/2} x\|}
{\sqrt{\E[\|\Sig^{-1/2} x\|^2]}}
= \frac{\|\Sig^{-1/2} x\|}{\sqrt{d}}
\leq \rho_0
.
\]
\end{remark}
\begin{remark}[Bounded covariates] \label{remark:leverage}
If $\|x\| \leq r$ almost surely, then
\[ \frac{\|(\Sig + \lambda I)^{-1/2}x\|}{\sqrt{\dl}}
\leq \frac{r}{\sqrt{(\inf\{\lambda_j\} + \lambda)\dl}}
\]
in which case Condition~\ref{cond:leverage} holds with $\rhol$ satisfying
\[
\rhol \leq \frac{r}{\sqrt{\lambda \dl}}
.
\]
\end{remark}

\subsubsection{Response model}

The response model considered in this work is a relaxation of the typical
Gaussian model; the model specifically allows for approximation error and
general subgaussian noise.
Define the random variables
\begin{equation} \label{eq:noise-bias}
\noise(x) := y - \E[y|x]
\quad \text{and} \quad
\bias(x) := \E[y|x] - \dotp{\beta,x}
\end{equation}
where $\noise(x)$ corresponds to the response noise, and $\bias(x)$
corresponds to the approximation error of $\beta$.
This gives the following modeling equation:
\[ y = \dotp{\beta,x} + \bias(x) + \noise(x) . \]
Conditioned on $x$, $\noise(x)$ is random, while $\bias(x)$ is
deterministic.

The noise is assumed to satisfy the following subgaussian moment condition:
\begin{condition}[Subgaussian noise] \label{cond:noise}
There exists finite $\snoise \geq 0$ such that, almost surely,
\[
\E\left[\exp(\eta \noise(x)) | x\right]
\leq \exp(\eta^2 \snoise^2/2)
\qquad \forall \eta \in \R
.
\]
\end{condition}
Condition~\ref{cond:noise} is satisfied, for instance, if $\noise(x)$ is
normally distributed with mean zero and variance $\snoise^2$.

For the next condition, define $\bl$ to be the minimizer of the
regularized mean squared error, \emph{i.e.},
\begin{equation} \label{eq:bl}
\bl := \arg\min_w \left\{ \E[(\dotp{w,x} - y)^2] + \lambda \|w\|^2
\right\} = (\Sig + \lambda I)^{-1} \E[xy] ,
\end{equation}
and also define
\begin{equation} \label{eq:biasl}
\biasl(x) := \E[y|x] - \dotp{\bl,x} .
\end{equation}
The final condition requires a bound on the size of $\biasl(x)$.
\begin{condition}[Bounded approximation error at $\lambda$]
\label{cond:bias}
There exist finite $\bbiasl \geq 0$ such that, almost surely,
\[
\frac{\|(\Sig + \lambda I)^{-1/2} x \biasl(x)\|}
{\sqrt{\E[\|(\Sig + \lambda I)^{-1/2} x\|^2]}}
= \frac{\|(\Sig + \lambda I)^{-1/2} x \biasl(x)\|}{\sqrt{\dl}}
\leq \bbiasl
.
\]
\end{condition}
The hard ``almost sure'' bound in Condition~\ref{cond:bias} can easily be
relaxed to moment conditions, but we do not consider it here for sake of
simplicity.
We also remark that $\bbiasl$ only appears in lower-order terms in the main
bounds.

\begin{remark}[Ordinary least squares]
\label{remark:bias-ols}
If $\lambda = 0$ and the dimension of the covariate space is $d$, then
Condition~\ref{cond:bias} reduces to the requirement that there exists a
finite $\bbias_0 \geq 0$ such that, almost surely,
\[
\frac{\|\Sig^{-1/2} x \bias(x)\|}
{\sqrt{\E[\|\Sig^{-1/2} x\|^2]}}
= \frac{\|\Sig^{-1/2} x \bias(x)\|}{\sqrt{d}}
\leq \bbias_0
.
\]
\end{remark}
\begin{remark}[Bounded approximation error] \label{remark:bias}
If $|\bias(x)| \leq a$ almost surely and Condition~\ref{cond:leverage}
(with parameter $\rhol$) holds, then
\begin{align*}
\frac{\|(\Sig + \lambda I)^{-1/2} x \biasl(x)\|}{\sqrt{\dl}}
& \leq \rhol |\biasl(x)| \\
& \leq \rhol (a + |\dotp{\beta-\bl,x}|) \\
& \leq \rhol (a + \|\beta-\bl\|_{\Sig+\lambda I} \|x\|_{(\Sig + \lambda
I)^{-1}}) \\
& \leq \rhol (a + \rhol \sqrt{\dl}\|\beta-\bl\|_{\Sig+\lambda I})
\end{align*}
where the first and last inequalities use Condition~\ref{cond:leverage},
the second inequality uses the definition of $\biasl(x)$
in~\eqref{eq:biasl} and the triangle inequality, and the third inequality
follows from Cauchy-Schwarz.
The quantity $\|\beta-\bl\|_{\Sig+\lambda I}$ can be bounded by
$\sqrt{\lambda} \|\beta\|$ using the arguments in the proof of
Proposition~\ref{proposition:regularization}.
In this case, Condition~\ref{cond:bias} is satisfied with
\[
\bbiasl \leq \rhol (a + \rhol \sqrt{\lambda\dl} \|\beta\|)
.
\]
If in addition $\|x\| \leq r$ almost surely, then
Condition~\ref{cond:leverage} and Condition~\ref{cond:bias} are satisfied
with
\[
\rhol \leq \frac{r}{\sqrt{\lambda \dl}}
\quad \text{and} \quad
\bbiasl \leq \rhol (a + r \|\beta\|)
\]
as per Remark~\ref{remark:leverage}.
\end{remark}

\subsection{Related work}
\label{section:previous-work}

The ridge and ordinary least squares estimators are classically studied in
the fixed design setting: the covariates $x_1, x_2, \dotsc, x_n$ are fixed
vectors in $\R^d$, and the responses $y_1, y_2, \dotsc, y_n$ are
independent random variables, each with mean $\E[y_i] = \dotp{\beta,x_i}$
and variance $\var(y_i) \leq \snoise^2$~\cite{lehmann-casella}.
The analysis reviewed in Section~\ref{section:fixed} reveals the
expected prediction error $\E[\|\blh - \beta\|_\Sig^2]$ is controlled by
the sum of a bias term, which is zero when $\lambda = 0$, and a variance
term, which is bounded by $\snoise^2\dll/n$.
As discussed in the introduction, our random design analysis of the ridge
estimator reveals the essential differences between fixed and random design
by comparing with this classical analysis.

Many classical analyses of the ridge and ordinary least squares estimators
in the random design setting (\emph{e.g.}, in the context of nonparametric
estimators) do not actually show nonasymptotic $O(d/n)$ convergence of the
mean squared error to that of the best linear predictor, where $d$ is the
dimension of the covariate space.
Rather, the error relative to the Bayes error is bounded by some multiple
$c > 1$ of the error of the optimal linear predictor relative to the Bayes
error, plus a $O(d/n)$ term~\cite{GKKW04}:
\[ \E[(\dotp{\bh,x}-\E[y|x])^2] \leq
c \cdot \E[(\dotp{\beta,x}-\E[y|x])^2] + O(d/n) . \]
Such bounds are appropriate in non-parametric settings where the error of
the optimal linear predictor also approaches the Bayes error at an $O(d/n)$
rate.
Beyond these classical results, analyses of ordinary least squares often
come with nonstandard restrictions on applicability or additional
dependencies on the spectrum of the second moment operator (see the recent
work of~\cite{AudCat10b} for a comprehensive survey of these results);
for instance, a result of~\cite{Catoni04} gives a bound
on the excess mean squared error of the form
\[ \|\bh - \beta\|_\Sig^2
\leq O\left(\frac{d + \log(\det(\hat\Sig)/\det(\Sig))}{n} \right) , \]
but the bound is only shown to hold when every linear predictor with low
empirical mean squared error satisfies certain boundedness conditions.

\ignore{
Some recent work on computational speed-ups for over-complete least squares
problems can be viewed as risk bounds for the ordinary least squares
estimator~\cite{DMMS10}.
However, because statistical risk bounds were not their objective, these
analyses do not separate noise and approximation error, have additional
$\log n$ factors in the bounds, and do not provide exponential tail
inequalities.
Our results for the ordinary least squares estimator can be readily applied
to their application of speeding-up least squares computations; we describe
this application in Section~\ref{section:applications}.
}

This work provides ridge regression bounds explicitly in terms of the
vector $\beta$ (as a sequence) and in terms of the eigenspectrum of the
second moment operator $\Sig$.
While the essential setting we study is not new, previous analyses make
unnecessarily strong boundedness assumptions or fail to give a bound in the
case $\lambda = 0$.
Here we review the analyses of \cite{Zhang05}, \cite{SmaZou07},
\cite{CD07}, and \cite{SHS09}.
\cite{Zhang05} assumes $\|x\| \leq b_x$ and
$|\dotp{\beta,x}-y| \leq b_{\bias}$ almost surely, and gives the bound
\[ \|\blh-\beta\|_\Sig^2 \leq \lambda
\|\blh-\beta\|^2 + c \cdot \frac{\dl \cdot (b_{\bias} + b_x
\|\blh-\beta\|)^2}{n} \]
for some $c>0$,
where $\dl$ is the effective dimension at scale $\lambda$ as defined
in~\eqref{eq:dim}.
The quantity $\|\blh-\beta\|$ is then bounded by assuming $\|\beta\| <
\infty$.
Thus, the dominant terms of the final bound have explicit dependences on
$b_{\bias}$ and $b_x$.
\cite{SmaZou07} assume that $|y| \leq b_y$ and $\|x\| \leq b_x$
almost surely, and prove the bound
\[ \|\blh - \bl\|_\Sig^2 \leq c' \cdot \frac{b_x^2b_y^2}{n\lambda^2} \]
for some $c'>0$ (and note that the bound becomes trivial when $\lambda =
0$); this is then used to bound $\|\blh-\beta\|_\Sig^2$ under explicit
assumptions on $\beta$.
\cite{CD07} assume $\|x\| \leq b_x$ almost surely, and prove the
bound (in their Theorem 4)
\[
  \|\blh-\beta\|_\Sig^2
  \leq c'' \cdot \biggl( \|\bl-\beta\|_\Sig^2
  + \frac{b_x \|\bl-\beta\|_\Sig^2}{n\lambda}
  + \frac{\snoise^2 \dl}{n}
  + o(1/n)
  \biggr)
  .
\]
Here, we also note that, if one desires the bound to hold with probability $\geq
1-e^{-t}$ for some $t>0$, then the leading factor $c''>1$ depends
quadratically on $t$.
Finally, \cite{SHS09} explicitly require $|y| \leq b_y$ and their
main bound on $\|\blh-\beta\|_\Sig^2$ (specialized for the ridge estimator)
depends on $b_y$ in a dominant term.
Moreover, this main bound contains $c''' \cdot ( \lambda\|\bl\|^2 +
\|\bl-\beta\|_\Sig^2 )$ as a dominant term for some $c'''>1$, and it is
only given under explicit decay conditions on the eigenspectrum (their
Equation 6).
The bound is also trivial when $\lambda = 0$.
Our result for ridge regression is given explicitly in terms of
$\|\bl-\beta\|_\Sig^2$ (and therefore explicitly in terms of $\beta$ as a
sequence, the eigenspectrum of $\Sig$, and $\lambda$); this quantity
vanishes when $\lambda=0$ and can be small even when $\|\beta\|$ itself is
large.
We note that $\|\bl-\beta\|_\Sig^2$ is precisely the bias term from the
classical fixed design analysis of ridge regression, and therefore is
natural to expect in a random design analysis.

Recently,~\cite{AudCat11} derived sharp risk bounds for the
ordinary least squares and ridge estimators (in addition to specially
developed PAC-Bayesian estimators) in a random design setting under very
mild moment assumptions using PAC-Bayesian techniques.
Their nonasymptotic bound for ordinary least squares holds with
probability at least $1-e^{-t}$ but only for $t \leq \ln n$; this is
essentially due to their weak moment assumptions.
By relying on stronger moment assumptions, we allow the probability tail
parameter $t$ to be as large as $\Omega(n/d)$.
Our analysis is also arguably more transparent and yields more reasonable
quantitative bounds.
The analysis of~\cite{AudCat11} for the ridge estimator is established
only in an asymptotic sense and therefore are not directly comparable to
those provided here.

Finally, although the focus of our present work is on understanding the
ordinary least squares and ridge estimators, it should also be mentioned
that a number of other estimators have been considered in the literature
with nonasymptotic prediction error bounds
\cite{Kolt06,AudCat11,HS13-heavy}.
Indeed, the works of \cite{AudCat11} and \cite{HS13-heavy} propose
estimators that require considerably weaker moment conditions on $x$ and
$y$ to obtain optimal rates.

\section{Random design regression}
\label{section:ridge}

This section presents the main results of the paper on the excess mean
squared error of the ridge estimator under random design (and its
specialization to the ordinary least squares estimator). First, we
review the standard fixed design analysis.

\subsection{Review of fixed design analysis}
\label{section:fixed}

It is informative to first review the fixed design analysis of the ridge
estimator.
Recall that, in this setting, the design points $x_1,x_2,\dotsc,x_n$ are
fixed (deterministic) vectors, and the responses $y_1,y_2,\dotsc,y_n$ are
independent random variables.
Therefore, we define $\Sig := \Sigh = n^{-1} \sum_{i=1}^n x_i \otimes x_i$
(which is nonrandom), and assume it has eigenvectors $\{v_j\}$ and
corresponding eigenvalues $\lambda_j := \dotp{v_j,\Sig v_j}$.
As in the random design setting, the linear function $\beta := \sum_j
\beta_j v_j$ where $\beta_j := (n \lambda_j)^{-1} \sum_{i=1}^n
\dotp{v_j,x_i} \E[y_i]$ minimizes the expected mean squared error,
\emph{i.e.},
\[ \beta := \arg\min_w \frac1n\sum_{i=1}^n \E[(\dotp{w,x_i}-y_i)^2] . \]
Similar to the random design setup, define $\noise(x_i) := y_i - \E[y_i]$
and $\bias(x_i) := \E[y_i] - \dotp{\beta,x_i}$ for $i=1,2,\dotsc,n$, so the
following modeling equation holds:
\[ y_i = \dotp{\beta,x_i} + \bias(x_i) + \noise(x_i) \]
for $i=1,2,\dotsc,n$.
Because $\Sig = \Sigh$, the ridge estimator $\blh$ in the fixed design
setting is an unbiased estimator of the minimizer of the regularized mean
squared error, \emph{i.e.},
\begin{align*}
\E[\blh]  
= (\Sig + \lambda I)^{-1} \left( \frac{1}{n} \sum_{i=1}^n
x_i\E[y_i] \right)
= \arg\min_w \left\{ \frac{1}{n} \sum_{i=1}^n \E[(\dotp{w,x_i} -
y_i)^2] + \lambda \|w\|^2 \right\}
.
\end{align*}
\ignore{
\begin{align*}
\E[\blh] & = (\Sig + \lambda I)^{-1} \left( \frac{1}{n} \sum_{i=1}^n
x_i\E[y_i] \right)
\\
& = \arg\min_w \left\{ \frac{1}{n} \sum_{i=1}^n \E[(\dotp{w,x_i} -
y_i)^2] + \lambda \|w\|^2 \right\}
.
\end{align*}}
This unbiasedness implies that the expected mean squared error of $\blh$
has the bias-variance decomposition
\begin{equation} \label{eq:bias-variance}
\E[\|\blh - \beta\|_\Sig^2]
= \|\E[\blh] - \beta\|_\Sig^2
+ \E[\|\blh - \E[\blh]\|_\Sig^2] .
\end{equation}
The following bound on the expected excess mean squared error easily
follows from this decomposition and the definition of $\beta$ (see,
\emph{e.g.}, Proposition~\ref{proposition:regularization}).
\begin{proposition}[Ridge regression: fixed design]
\label{proposition:fixed}
Fix $\lambda \geq 0$, and assume $\Sig + \lambda I$ is invertible.
If there exists $\snoise \geq 0$ such that $\var(y_i^2) \leq \snoise^2$ for
all $i=1,2,\dotsc,n$, then
\[
\E[\|\blh - \beta\|_\Sig^2]
\leq \sum_j \frac{\lambda_j}{(\frac{\lambda_j}{\lambda} + 1)^2} \beta_j^2
+ \frac{\snoise^2}{n}
\sum_j \left( \frac{\lambda_j}{\lambda_j + \lambda} \right)^2
\]
with equality iff $\var(y_i) = \snoise^2$ for all $i=1,2,\dotsc,n$.
\end{proposition}
\begin{remark}[Effect of approximation error in fixed design]
\label{remark:fixed-approx}
Observe that $\bias(x_i)$ has no effect on the expected excess mean squared
error.
\end{remark}
\begin{remark}[Effective dimension]
\label{remark:fixed-dimension}
The second sum in the bound is equal to $\dll$, a notion of effective
dimension at regularization level $\lambda$.
\end{remark}
\begin{remark}[Ordinary least squares in fixed design]
\label{remark:fixed-ols}
Setting $\lambda = 0$ gives the
following bound for the ordinary least squares estimator $\bh_0$:
\[ \E[\|\bh_0 - \beta\|_\Sig^2] \leq \frac{\snoise^2 d}{n} \]
where, as before, equality holds iff $\var(y_i) = \snoise^2$ for all $i =
1,2,\dotsc,n$.
\end{remark}

\subsection{Ordinary least squares}

Our analysis of the ordinary least squares estimator (under random
design) is based on a simple decomposition of the excess mean squared
error, similar to the one from the fixed design analysis.  To state the
decomposition, first let $\bb_0$ denote the conditional expectation of
the least squares estimator $\bh_0$ conditioned on $x_1,x_2,\dotsc,x_n$, \emph{i.e.},
\begin{equation*} 
\bb_0 := \E[\bh_0 | x_1,x_2,\dotsc,x_n]
= \Sigh^{-1} \Eh[x\E[y|x]]
.
\end{equation*}
Also, define the bias and variance as:
\[
\epsII := \|\bb_0 - \beta\|_\Sig^2  \ ,
\quad \quad
\epsIII := \|\bh_0 - \bb_0\|_\Sig^2 
\]

\begin{proposition}[Random design decomposition]
\label{proposition:ols_decomp}
We have
\begin{align*}
\|\bols - \beta\|_\Sig^2
& \leq \epsII + 2\sqrt{\epsII\epsIII} + \epsIII
\\
& \leq 2(\epsII + \epsIII)
\end{align*}
\ignore{
where
$\bols$ is the ordinary least squares estimator (see
~\eqref{eq:blh} with $\lambda = 0$) and
$\beta$ achieves the minimal least squared error (see~\eqref{eq:beta}).
}
\end{proposition}
\begin{proof}
The claim follows from the triangle inequality and the fact $(a+b)^2 \leq
2(a^2+b^2)$.
\end{proof}

\begin{remark}
Note that, in general, $\E[\bh_0] \neq \beta$ (unlike in the fixed
design setting where $\E[\bh_0] =\beta$). Hence, our
decomposition differs from that in the fixed design analysis
(see~\eqref{eq:bias-variance}).
\end{remark}

\ignore{
Theorem~\ref{theorem:ridge} readily applies the ordinary least squares
estimator simply by setting $\lambda = 0$.
This requires the dimension of the covariate space $d$ to be finite.
However, a slightly different but more direct argument can be made in this
case, and results in a simpler bound with fewer lower-order terms.
}

Our first main result characterizes the excess loss of the ordinary
least squares estimator.

\begin{theorem}[Ordinary least squares regression] \label{theorem:ols}
Pick any $t > \max\{0, 2.6 - \log d \}$.
Assume 
Condition~\ref{cond:leverage} (with
parameter $\rho_0$), Condition~\ref{cond:noise} (with $\snoise$),
and Condition~\ref{cond:bias} (with $\bbias_0$) hold and that
\[ n \geq 6\rho_0^2 d(\log d + t) . \]
With probability at least $1 - 3e^{-t}$, the following holds:
\begin{enumerate}
\item \underline{Relative spectral norm error in $\Sigh$}:
$\Sigh$ is invertible, and
\[ \|\Sig^{1/2}\Sigh^{-1}\Sig^{1/2}\| \leq (1 - \deltaI)^{-1} , \]
where
$\Sig$ is defined in~\eqref{eq:Sig},
$\Sigh$ is defined in~\eqref{eq:Sigh},
and
\[
\deltaI := 
\sqrt{\frac{4\rho_0^2d(\log d + t)}{n}} + \frac{2\rho_0^2d(\log d + t)}{3n}
\]
(note that the lower bound on $n$ ensures $\deltaI \leq 0.93 < 1$).

\item \underline{Effect of bias due to random design}:
\begin{align*}
\epsII & \leq
\frac{2}{(1-\deltaI)^2}
\Biggl(
\frac{\E[\|\Sig^{-1/2}x\bias(x)\|^2]}{n}
(1 + \sqrt{8t})^2
+ \frac{16\bbias_0^2 dt^2}{9n^2}
\Biggr)
\\
& \hphantom{:}\leq
\frac{2}{(1-\deltaI)^2}
\Biggl(
\frac{\rho_0^2 d \E[\bias(x)^2]}{n} (1 + \sqrt{8t})^2
+ \frac{16\bbias_0^2 dt^2}{9n^2}
\Biggr)
,
\end{align*}
and $\bias(x)$ is defined in~\eqref{eq:biasl}.

\item \underline{Effect of noise}:
\ignore{
\[ \|\bh_0 - \bb_0\|_\Sig^2 \leq \epsIII \]
where
}
\[
\epsIII \leq
\frac1{1-\deltaI}
\cdot
\frac{\snoise^2 (d + 2 \sqrt{d t} + 2 t)}{n}
.
\]

\end{enumerate}
\end{theorem}

\begin{remark}[Simplified form] \label{remark:simplified-ols}
Suppressing the terms that are $o(1/n)$, the overall bound from
Theorem~\ref{theorem:ols} is
\[
\|\bh_0 - \beta\|_\Sig^2
\leq
\frac{2\E[\|\Sig^{-1/2}x\bias(x)\|^2]}{n}
(1 + \sqrt{8t})^2 
+
\frac{\snoise^2 (d + 2 \sqrt{d t} + 2 t)}{n}
+ o(1/n)
\]
(so $\bbias_0$ appears only in the $o(1/n)$ terms).
If the linear model is correct (\emph{i.e.}, $\E[y|x] = \dotp{\beta,x}$
almost surely), then
\begin{equation} \label{eq:simple-ols}
\|\bh_0 - \beta\|_\Sig^2
\leq \frac{\snoise^2 (d + 2 \sqrt{d t} + 2 t)}{n}
+ o(1/n)
.
\end{equation}
One can show that the constants in the first-order term in~\eqref{eq:simple-ols}
are the same as those that one would obtain for a fixed design tail bound.
\end{remark}

\begin{remark}[Tightness of the bound] \label{remark:tightness-ols}
Since
\[
\|\bb_0 - \beta\|_\Sig^2
= \|(\Sig^{1/2} \Sigh^{-1} \Sig^{1/2}) \Eh[\Sig^{-1/2} x \bias(x)]\|^2
\]
and
\[ \|\Sig^{1/2} \Sigh^{-1} \Sig^{1/2} - I\| \to 0 \]
as $n \to \infty$ (Lemma~\ref{lemma:deltal-spectral}), $\|\bb_0 -
\beta\|_\Sig^2$ is within constant factors of $\|\Eh[\Sig^{-1/2}
x\bias(x)]\|^2$ for sufficiently large $n$.
Moreover,
\[ \E[\|\Eh[\Sig^{-1/2} x\bias(x)]\|^2] = \frac{\E[\|\Sig^{-1/2}
x\bias(x)\|^2]}{n} , \]
which is the main term that appears in the bound for $\epsII$.
Similarly, $\|\bh_0 - \bb_0\|_\Sig^2$ is within constant factors of
$\|\bh_0 - \bb_0\|_{\Sigh}^2$ for sufficiently large $n$, and
\[ \E[\|\bh_0 - \bb_0\|_{\Sigh}^2] \leq \frac{\snoise^2 d}{n} \]
with equality iff $\var(y) = \snoise^2$ (this comes from the fixed design
risk bound in Remark~\ref{remark:fixed-ols}).
Therefore, in this case where $\var(y) = \snoise^2$, we conclude that the
bound Theorem~\ref{theorem:ols} is tight up to constant factors and
lower-order terms.
\end{remark}

\subsection{Random design ridge regression}

The analysis of the ridge estimator under random design is again based on a
simple decomposition of the excess mean squared error.
Here, let $\blb$ denote the conditional
expectation of $\blh$ given $x_1,x_2,\dotsc,x_n$, \emph{i.e.},
\begin{equation} \label{eq:blb}
\blb := \E[\blh | x_1,x_2,\dotsc,x_n]
= (\Sigh + \lambda I)^{-1} \Eh[x\E[y|x]]
.
\end{equation}
Define the bias from regularization, the bias from the random
design, and the variance as:
\[
\epsI := \|\bl - \beta\|_\Sig^2 , \quad \quad 
\epsII := \|\blb - \bl\|_\Sig^2 , \quad \quad 
\epsIII := \ \|\blh - \blb\|_\Sig^2 ,
\] 
where $\bl$ is the minimizer of the
regularized mean squared error (see~\eqref{eq:bl}).

\begin{proposition}[General random design decomposition]
\label{proposition:decomposition}
\begin{align*}
\|\blh - \beta\|_\Sig^2
& \leq \epsI + \epsII + \epsIII
+ 2(\sqrt{\epsI\epsII} + \sqrt{\epsI\epsIII} + \sqrt{\epsII\epsIII})
\\
& \leq 3(\epsI + \epsII + \epsIII)
\end{align*}
\ignore{
\begin{align*}
\|\blh - \beta\|_\Sig^2
& \leq
\|\bl - \beta\|_\Sig^2
+ \|\blb - \bl\|_\Sig^2
+ \|\blh - \blb\|_\Sig^2
\\
& \quad{}
+ 2 \|\bl - \beta\|_\Sig \|\blb - \bl\|_\Sig
+ 2 \|\bl - \beta\|_\Sig \|\blh - \blb\|_\Sig
+ 2 \|\blb - \bl\|_\Sig \|\blh - \blb\|_\Sig
\\
& \leq
3 \left( \|\bl - \beta\|_\Sig^2
+ \|\blb - \bl\|_\Sig^2
+ \|\blh - \blb\|_\Sig^2 \right)
.
\end{align*}
}
\end{proposition}
\begin{proof}
The claim follows from the triangle inequality and the fact $(a+b)^2 \leq
2(a^2+b^2)$.
\end{proof}

\begin{remark}
Again, note that $\E[\blh] \neq \bl$ in general, so the bias-variance
decomposition in~\eqref{eq:bias-variance} from the fixed design analysis is
not directly applicable in the random design setting.
\end{remark}

The following theorem is the main result of the paper:
\begin{theorem}[Ridge regression] \label{theorem:ridge}
Fix some $\lambda \geq 0$, and pick any $t > \max\{0, 2.6 - \log \dlt \}$.
Assume Condition~\ref{cond:leverage} (with
parameter $\rhol$), Condition~\ref{cond:noise} (with parameter $\snoise$),
and Condition~\ref{cond:bias} (with parameter $\bbiasl$) hold; and that
\[ n \geq 6\rhol^2\dl(\log\dlt + t) , \]
where $d_{p,\lambda}$ for $p \in \{1,2\}$ is defined in~\eqref{eq:dim}, and
$\dlt$ is defined in~\eqref{eq:dlt}.

With probability at least $1 - 4e^{-t}$, the following holds:
\begin{enumerate}
\ignore{
\item \underline{Overall bound on excess mean squared error}:
\begin{align*}
\|\blh - \beta\|_\Sig^2
& \leq \epsI + \epsII + \epsIII
+ 2(\sqrt{\epsI\epsII} + \sqrt{\epsI\epsIII} + \sqrt{\epsII\epsIII})
\\
& \leq 3(\epsI + \epsII + \epsIII)
\end{align*}
where
$\blh$ is the ridge estimator defined in~\eqref{eq:blh};
$\beta$ is the least squares linear predictor defined in~\eqref{eq:beta};
and $\epsI$, $\epsII$, and $\epsIII$ are given below.
}
\item \underline{Relative spectral norm error in $\Sigh + \lambda I$}:
$\Sigh + \lambda I$ is invertible, and
\[ \|(\Sig + \lambda I)^{1/2}(\Sigh + \lambda I)^{-1}(\Sig + \lambda
I)^{1/2}\| \leq (1 - \deltaI)^{-1} , \]
where
$\Sig$ is defined in~\eqref{eq:Sig},
$\Sigh$ is defined in~\eqref{eq:Sigh},
and
\[
\deltaI := 
\sqrt{\frac{4\rhol^2\dl(\log\dlt + t)}{n}} + \frac{2\rhol^2\dl (\log\dlt +
t)}{3n}
\]
(note that the lower bound on $n$ ensures $\deltaI \leq 0.93 < 1$).

\item \underline{Frobenius norm error in $\Sigh$}:
\[ \|(\Sig + \lambda I)^{-1/2}(\Sigh - \Sig)(\Sig + \lambda
I)^{-1/2}\|_\F \leq \sqrt{\dl} \deltaII , \]
where
\[
\deltaII := 
\sqrt{\frac{\rhol^2\dl - \dll/\dl}{n}}(1+\sqrt{8t})
+ \frac{4\sqrt{\rhol^4\dl + \dll/\dl}t}{3n}
.
\]

\item \underline{Effect of regularization}:
\ignore{
\[ \|\bl - \beta\|_\Sig^2 = \epsI \]
where $\bl$ is defined in~\eqref{eq:bl}, and
}
\[
\epsI \leq \sum_j \frac{\lambda_j}{(\frac{\lambda_j}{\lambda} + 1)^2} \beta_j^2
. \]
If $\lambda = 0$, then $\epsI = 0$.

\item \underline{Effect of bias due to random design}:
\ignore{
\[ \|\blb - \bl\|_\Sig^2 \leq \epsII \]
where $\blb$ is defined in~\eqref{eq:blb},
}
\begin{align*}
\epsII & \leq
\frac{2}{(1-\deltaI)^2}
\Biggl(
\frac{\E[\|(\Sig + \lambda I)^{-1/2}(x\biasl(x) - \lambda\bl)\|^2]}{n}
(1 + \sqrt{8t})^2
+ \frac{16\bigl(\bbiasl \sqrt{\dl} + \sqrt{\epsI} \bigr)^2t^2}{9n^2}
\Biggr)
\\
& \hphantom{:}\leq
\frac{4}{(1-\deltaI)^2}
\Biggl(
\frac{\rhol^2 \dl \E[\biasl(x)^2] + \epsI}{n} (1 + \sqrt{8t})^2
+ \frac{\bigl(\bbiasl \sqrt{\dl} + \sqrt{\epsI} \bigr)^2t^2}{n^2}
\Biggr)
,
\end{align*}
and $\biasl(x)$ is defined in~\eqref{eq:biasl}.
If $\lambda = 0$, then $\biasl(x) = \bias(x)$ as defined
in~\eqref{eq:noise-bias}.

\item \underline{Effect of noise}:
\ignore{
\[ \|\blh - \blb\|_\Sig^2 \leq \epsIII \]
where
}
\[
\epsIII \leq
\frac{\snoise^2 \Bigl(\dll + \sqrt{\dl\dll}\deltaII
\Bigr)}{n (1-\deltaI)^2}
+ \frac{2\snoise^2 \sqrt{\Bigl(\dll + \sqrt{\dl\dll}\deltaII\Bigr)t}}{n
(1-\deltaI)^{3/2}}
+ \frac{2\snoise^2 t}{n (1-\deltaI)}
.
\]

\end{enumerate}
\end{theorem}

We now discuss various aspects of Theorem~\ref{theorem:ridge}.
\begin{remark}[Simplified form] \label{remark:simplified}
Ignoring the terms that are $o(1/n)$ and treating $t$ as a constant, the
overall bound from Theorem~\ref{theorem:ridge} is
\begin{align*}
\|\blh - \beta\|_\Sig^2
& \leq \|\bl-\beta\|_\Sig^2
+ O\left( \frac{\E[\|(\Sig + \lambda I)^{-1/2}(x\biasl(x) -
\lambda\bl)\|^2] + \snoise^2 \dll}{n} \right)
\\
& \leq \|\bl-\beta\|_\Sig^2
+ O\left( \frac{\rhol^2\dl\E[\biasl(x)^2] + \|\bl-\beta\|_\Sig^2
+ \snoise^2 \dll}{n} \right)
\\
& \leq \|\bl-\beta\|_\Sig^2
+ O\left( \frac{\rhol^2\dl\E[\bias(x)^2]
+ (\rhol^2\dl+1)\|\bl-\beta\|_\Sig^2
+ \snoise^2 \dll}{n} \right)
\end{align*}
where the last inequality follows from the fact $\sqrt{\E[\biasl(x)^2]}
\leq \sqrt{\E[\bias(x)^2]} + \|\bl-\beta\|_\Sig$.
\end{remark}
\begin{remark}[Effect of errors in $\Sigh$] \label{remark:Sigma}
The accuracy of $\Sigh$ has a relatively mild effect on the bound---it
appears essentially through multiplicative factors $(1-\deltaI)^{-1} = 1 +
O(\deltaI)$ and $1+\deltaII$, where both $\deltaI$ and $\deltaII$ are
decreasing with $n$ (as $n^{-1/2}$), and therefore only contribute to
lower-order terms overall.
\end{remark}
\begin{remark}[Effect of approximation error] \label{remark:approx}
The effect of approximation error is isolated in the term
$\|\blb-\bl\|_\Sig^2$.
The bound $\epsII$ scales with a fourth-moment quantity $\E[\|(\Sig +
\lambda I)^{-1/2}(x\biasl(x) - \lambda\bl\|^2]$; when using the looser
bound $O(\rhol^2\dl\E[\bias(x)^2] + (\rho^2\dl+1) \|\bl-\beta\|_\Sig^2)$,
the overall simplified bound from Remark~\ref{remark:simplified} can be
viewed as
\begin{align*}
  \lefteqn{
\E[(\dotp{\blh,x}-\E[y|x])^2 | \blh]
} \\
& \leq \E[(\dotp{\beta,x}-\E[y|x])^2]
\left( 1 + \frac{c_1 \rhol^2\dl}{n} \right)
+ \E[\dotp{\bl-\beta,x}^2] \left( 1 + \frac{c_2 (\rhol^2\dl+1)}{n} \right)
\\
& \quad{}
+ \text{terms due to stochastic noise}
\end{align*}
for some positive constants $c_1$ and $c_2$.
Therefore, the (bound on the) mean squared error of $\blh$ is the sum of
two contributions (up to lower-order terms): the first is a scaling of the
approximation errors $\E[(\dotp{\beta,x} - \E[y|x])^2] +
\E[\dotp{\bl-\beta,x}^2]$, where the scaling $1 + O((\rhol^2\dl+1)/n)$
tends to one as $n \to \infty$; and the second is the stochastic noise
contribution.
The approximation error contribution is unique to random design, while the
stochastic noise appears in both random and fixed design.
\end{remark}
\begin{remark}[Bounded covariates] \label{remark:bounded}
  Suppose $\bias(x) = 0$ and that there exists $r>0$ such that $\|x\| \leq
  r$ almost surely.
  This is the setting of a well-specified model with bounded covariates;
  the minimax risk over the class of models $\beta$ with $\|\beta\| \leq B$
  for some $B>0$ is at least $\Omega(\sqrt{\snoise^2 r^2B^2 /
  n})$~\cite{Nussbaum99}.
  In this case, using the inequalities $\|\bl-\beta\|_\Sig^2 \leq
  \lambda\|\beta\|^2/2$ and $\dll \leq \tr(\Sig)/(2\lambda)$, the
  simplified bound from Remark~\ref{remark:simplified} reduces to
  \[
    \|\blh - \beta\|_\Sig^2 \leq
    \Biggl( 1 + O\biggl( \frac{1 + r^2/\lambda}{n} \biggr) \Biggr)
    \cdot \frac{\lambda\|\beta\|^2}{2} + \frac{\snoise^2}{n} \cdot
    \frac{\tr(\Sig)}{2\lambda} .
  \]
  Choosing $\lambda > 0$ to minimize the bound and using the fact
  $\tr(\Sig) \leq r^2$ gives
  \[
    \|\blh - \beta\|_\Sig^2
    \leq \sqrt{\frac{\snoise^2r^2B^2}{n} \cdot \biggl( 1 + O(1/n) \biggr)}
    + O\biggl( \frac{r^2B^2}{n} \biggr) ,
  \]
  which matches the lower bound up to constant factors and lower-order
  terms.
\end{remark}

\begin{remark}[Application to smoothing splines] \label{section:splines}
The applications of ridge regression considered by~\cite{Zhang05} can also
be analyzed using Theorem~\ref{theorem:ridge} (although technically our
result is only proved in the finite-dimensional setting).
We specifically consider the problem of approximating a periodic function
with smoothing splines, which are functions $f \colon \R \to \R$ whose
$s$-th derivatives $f^{(s)}$, for some $s > 1/2$, satisfy
\[ \int \left( f^{(s)}(t) \right)^2 dt < \infty . \]
The one-dimensional covariate $t \in \R$ can be mapped to the infinite-dimensional representation $x := \phi(t) \in \R^\infty$ where
\[
x_{2k} := \frac{\sin(kt)}{(k+1)^s}
\quad\text{and}\quad
x_{2k+1} := \frac{\cos(kt)}{(k+1)^s}
,
\quad k \in \{0, 1, 2, \dotsc \}
.
\]
Assume that the regression function is
\[ \E[y|x] = \dotp{\beta,x} \]
so $\bias(x) = 0$ almost surely.
Observe that $\|x\|^2 \leq \frac{2s}{2s-1}$, so
Condition~\ref{cond:leverage} is satisfied with
\[
\rhol := \left(\frac{2s}{2s-1}\right)^{1/2} \frac1{\sqrt{\lambda\dl}}
\]
as per Remark~\ref{remark:leverage}.
Therefore, the simplified bound from Remark~\ref{remark:simplified} becomes
in this case
\begin{align*}
\|\blh - \beta\|_\Sig^2
& \leq \|\bl - \beta\|_\Sig^2
+ C \cdot \left( \frac{2s}{2s-1} \cdot \frac{
\|\bl-\beta\|_\Sig^2}{\lambda n}
+ \frac{\|\bl - \beta\|_\Sig^2 + \snoise^2\dll}{n} \right)
\\
& \leq \frac{\lambda\|\beta\|^2}{2}
+ C \cdot \frac{\snoise^2\dll}{n}
+ C \cdot \left( \frac{2s}{2s-1} + \frac{\lambda}{2} \right) \cdot
\frac{\|\beta\|^2}{n}
\end{align*}
for some constant $C > 0$, where we have used the inequality $\|\bl -
\beta\|_\Sig^2 \leq \lambda \|\beta\|^2 / 2$.
\cite{Zhang05} shows that
\[ \dl \leq \inf_{k \geq 1} \left\{ 2k + \frac{2/\lambda}{(2s-1)k^{2s-1}}
\right\} . \]
Since $\dll \leq \dl$, it follows that setting $\lambda := k^{-2s}$ where
$k = \lfloor  ( (2s-1)n / (2s) )^{1/(2s+1)} \rfloor$
gives the bound
\[
\|\blh - \beta\|_\Sig^2
\leq
\left( \frac{\|\beta\|^2}{2} + 2C\snoise^2 \right)
\cdot \left( \frac{2s-1}{2s} \cdot n
\right)^{-\frac{2s}{2s+1}}
+ \text{lower-order terms}
\]
which has the optimal data-dependent rate of
$n^{-\frac{2s}{2s+1}}$~\cite{Stone82}.
\end{remark}
\begin{remark}[Comparison with fixed design] \label{remark:fixed}
As already discussed, the ridge estimator behaves similarly under fixed and
random designs, with the main differences being the lack of errors in
$\Sigh$ under fixed design, and the influence of approximation error under
random design.
These are revealed through the quantities $\rhol$ and $\dl$ (and $\bbiasl$
in lower-order terms), which are needed to apply the probability tail
inequalities.
Therefore, the scaling of $\rhol^2\dl$ with $\lambda$ crucially controls
the effect of random design compared with fixed design.
\end{remark}

\section{Application to accelerating least squares computations}
\label{section:applications}

Our results for the ordinary least squares estimator can be used to analyze
a randomized approximation scheme for overcomplete least squares
problems~\cite{DMMS10,RokTyg08}.
The goal of these randomized methods is to approximately solve the least
squares problem
\[ \min_{w \in \R^d} \frac1m \|Aw - b\|^2 \]
for some large, full-rank design matrix $A \in \R^{m \times d}$ ($m \gg d$)
and vector $b \in \R^m$.
Note that using a standard method to exactly solve the least squares
problem requires $\Omega(md^2)$ operations, which can be prohibitive for
large-scale problems.
However, when an approximate solution is satisfactory, significant
computational savings can be achieved through the use of randomization.

\subsection{A randomized approximation scheme for least squares}

The approximation scheme is as follows:
\begin{enumerate}
\item The columns of $A$ and the vector $b$ are first subjected to a
randomly chosen rotation matrix (\emph{i.e.}, an orthogonal transformation)
$\Theta \in \R^{m \times m}$.
The distribution over rotation matrices that may be used is discussed
below.

\item A sample of $n$ rows of $[\Theta A, \Theta b] \in \R^{m \times
(d+1)}$ are then selected uniformly at random with replacement;
let $\{ [x_i^\t, y_i] : i = 1,2,\dotsc,n \}$ (where $x_i \in \R^d$ and $y_i
\in \R$) be the $n$ selected rows of $[\Theta A, \Theta b]$.

\item Finally, the least squares problem
\[ \min_{w \in \R^d} \frac1n \sum_{i=1}^n (\dotp{w,x_i} -
y_i)^2 \]
is solved by computing the ordinary least squares estimator $\bh_0$ on the
sample $\{ (x_i,y_i) : i = 1,2,\dotsc,n \}$.

\end{enumerate}
The motivation for the random rotation $\Theta$ is captured in
Lemma~\ref{lemma:2->infty}, which shows that, if $\Theta$ is chosen randomly from
certain distributions over rotation matrices, then applying $\Theta$ to $A$
and $b$ creates an equivalent least squares problem for which the
statistical leverage parameter (the quantity $\rho_0$ in
Condition~\ref{cond:leverage}) is small.
Consequently, the new least squares problem can be approximately solved
with a small random sample, as per Theorems~\ref{theorem:ridge} and
\ref{theorem:ols}.
Without the random rotation, the statistical leverage parameter could be so
large that small random sample of the rows will likely miss a row crucial
for obtaining an accurate approximation.
The role of statistical leverage in this setting was also pointed out
by~\cite{DM10}, although Lemma~\ref{lemma:2->infty} makes the connection
more direct.
We note that Lemma~\ref{lemma:2->infty} and the analysis below can be
generalized to the case where $\Theta$ is only approximately orthogonal;
for most standard distributions over rotation matrices, the additional
error terms that arise do not affect the overall analysis.

The running time of the approximation scheme is given by (i) the time
required to apply the $m \times m$ random rotation operator $\Theta$ to the
original $m \times (d+1)$ matrix $[A, b]$ and randomly sample $n$ rows,
plus (ii) the time to solve the least squares problem on the smaller design
matrix of size $n \times d$.
For (i), na\"ively applying an arbitrary $m \times m$ rotation matrix
requires $\Omega(m^2 d)$ operations; however, there are (distributions
over) rotation matrices for which this running time can be reduced to $O(md
\log m)$ (see Example~\ref{example:hadamard} in
Section~\ref{section:random-rotations} below), which is a considerable
speed-up when $m$ is large.
In fact, because only $n$ out of $m$ rows are to be retained anyway, this
computation can be reduced to $O(md \log n)$~\cite{AilCha09}.
For (ii), standard methods can produce the ordinary least squares estimator
or the ridge regression estimator with $O(nd^2)$ operations.
Therefore, we are interested in the sample size $n$ that suffices to yield
an accurate approximation.

\subsection{Analysis of the approximation scheme}

Our approach to analyzing the above approximation scheme is to treat it as
a random design regression problem.
We apply Theorem~\ref{theorem:ols} in this
setting to give error bounds for the solution produced by the approximation
scheme.

Let $(x,y) \in \R^d \times \R$ be a random pair distributed uniformly over
the rows of $[\Theta A, \Theta b]$, where we assume that $\Theta$ is
randomly chosen from a suitable distribution over rotation matrices such as
those described in Example~\ref{example:uniform} and
Example~\ref{example:hadamard}.
Lemma~\ref{lemma:2->infty} (below) implies that
there exists a constant $c_0 > 0$ such that Condition~\ref{cond:leverage}
is satisfied with
\[ \rho_0^2 \leq c_0 \cdot \left(1 + \frac{\log m + \tau}{d} \right) \]
with probability at least $1-e^{-\tau}$ over the choice of the random
rotation matrix $\Theta$.
Henceforth, we condition on the event that this holds.

Let $\beta \in \R^d$ be the solution to the original least squares problem
(\emph{i.e.}, $\beta := \arg\min_w \|Aw-b\|^2 / m$), and let $\bh_0 \in
\R^d$ be the ordinary least squares estimator computed on the random sample of the rows of
$[\Theta A, \Theta b]$.
Note that, for any $w \in \R^d$,
\[
\E[(\dotp{w,x}-y)^2]
= \frac1m \|\Theta Aw-\Theta b\|^2
= \frac1m \|Aw-b\|^2
.
\]
Moreover, we may assume for simplicity that $y - \dotp{\beta,x} = \bias(x)$
(\emph{i.e.}, there is no stochastic noise), so $\E[\bias(x)^2] =
\E[(\dotp{\beta,x}-y)^2] = \|A\beta-b\|^2 / m$.

By Theorem~\ref{theorem:ols}, if at least
\[ n \geq
6 \bigl(d + c_0 (\log m + \tau) \bigr) (\log d + t)
\]
rows of $[\Theta A,\Theta b]$ are sampled, then the ordinary least
squares estimator $\bh_0$ satisfies the following approximation error
guarantee (with probability at least $1-3e^{-t}$ over the random sample of
rows):
\[
\frac1m \|A\bh_0 - b\|^2
\leq \frac1m \|A\beta - b\|^2
\cdot \left( 1 + c_1 \frac{(d + \log m + \tau) t}{n} \right) +
o(1/n)
\]
for some constant $c_1 > 0$.
We note that the $o(1/n)$ terms can be removed if one only requires
constant probability of success (\emph{i.e.}, $\tau$ and $t$ treated as
constants), as is considered by~\cite{DMMS10}.
In this case, we achieve an error bound of
\[
\frac1m \|A\bh_0 - b\|^2
\leq \frac1m \|A\beta - b\|^2 \cdot (1 + \epsilon)
\]
for $\epsilon > 0$ provided that the number of rows sampled is
\[ n \geq c_2 (d + \log m) \left( \frac1{\epsilon} + \log d
\right) \]
for some constant $c_2 > 0$.

\subsection{Random rotations and bounding statistical leverage}
\label{section:random-rotations}

The following lemma gives a simple condition on the distribution of the
random orthogonal matrix $\Theta \in \R^{n \times n}$ used to preprocess a
data matrix $A$ so that Condition~\ref{cond:leverage} is applicable to a
random vector $x$ drawn uniformly from the rows of $\Theta A$.
Its proof is a straightforward application of
Lemma~\ref{lemma:quadratic}.

\begin{lemma} \label{lemma:2->infty}
Fix any $\tau > 0$ and $\lambda \geq 0$.
Suppose $\Theta \in \R^{m \times m}$ is a random orthogonal matrix
and $\kappa > 0$ is a constant such that
\begin{equation} \label{eq:rotation-cond}
\E\left[ \exp\left( \alpha^\t \bigl(\sqrt{m} \Theta^\t e_i\bigr)
\right) \right] \leq \exp\left(\kappa \|\alpha\|^2/2 \right)
, \quad
\forall \alpha \in \R^m ,
\forall i=1,2,\dotsc,m ,
\end{equation}
where $e_i$ is the $i$-th coordinate vector in $\R^m$.
Let $A \in \R^{m \times d}$ be any matrix of rank $d$, and let $\Sig :=
(1/m) (\Theta A)^\t (\Theta A) = (1/m) A^\t A$.
There exists
\[
\rhol^2 \leq \kappa \left( 1 + 2\sqrt{\frac{\log m + \tau}{\dl}} +
\frac{2(\log m + \tau)}{\dl} \right)
\]
such that
\[
\Pr\left[ \max_{i=1,2,\dotsc,m} \|(\Sig + \lambda I)^{-1/2} (\Theta A)^\t
e_i \|^2 > \rhol^2 \dl \right] \leq e^{-\tau}
\]
where
$\dl := \sum_{j=1}^d \frac{\lambda_j}{\lambda_j + \lambda}$
and $\{ \lambda_1, \lambda_2, \dotsc, \lambda_d \}$ are the eigenvalues of
$\Sig$.
\end{lemma}
\begin{proof}
Let $z_i := \sqrt{m} \Theta^\t e_i$ for each $i=1,2,\dotsc,n$.
Let $U \in \R^{m \times d}$ be a matrix of left orthonormal singular
vectors of $(1/\sqrt{m}) A$, and let $D_\lambda :=
\diag(\frac{\lambda_1}{\lambda_1 + \lambda}, \frac{\lambda_2}{\lambda_2 +
\lambda}, \dotsc, \frac{\lambda_d}{\lambda_d + \lambda})$.
Note that $D_\lambda = I$ if $\lambda = 0$.
We have
\[
\|(\Sig + \lambda I)^{-1/2} (\Theta A)^\t e_i\|
= \|\sqrt{m} D_\lambda^{1/2} U^\t \Theta^\t e_i\|
= \|D_\lambda^{1/2} U^\t z_i\|
.
\]
Since $\tr(UD_\lambda U^\top) = \dl$, $\tr(UD_\lambda^2 U^\top) \leq \dl$,
and $\lambda_{\max}[UD_\lambda U^\top] \leq 1$, Lemma~\ref{lemma:quadratic}
implies
\[
\Pr\left[ \|D_\lambda^{1/2} U^\t z_i\|^2 > \kappa \left( \dl + 2\sqrt{\dl
(\log m + \tau)} + 2(\log m + \tau) \right) \right] \leq e^{-\tau}/m
.
\]
Therefore, by a union bound,
\[
\Pr\left[
\max_{i=1,2,\dotsc,m} \|(\Sig + \lambda I)^{-1/2} (\Theta A)^\t e_i\|^2
> \kappa \left( \dl + 2\sqrt{\dl (\log m + \tau)} +
2(\log m + \tau) \right) \right] \leq e^{-\tau}
.
\qedhere
\]
\end{proof}

Below, we give two simple examples under which the
condition~\eqref{eq:rotation-cond} in Lemma~\ref{lemma:2->infty} holds.
\begin{example} \label{example:uniform}
Let $\Theta$ be distributed uniformly over all $m \times m$ orthogonal
matrices.
Fix any $i=1,2,\dotsc,m$.
The random vector $v := \Theta^\t e_i$ is distributed uniformly on the
unit sphere $\S^{m-1}$.
Let $l$ be a $\chi$ random variable with $m$ degrees of freedom, so $z :=
lv$ follows an isotropic multivariate Gaussian distribution.
By Jensen's inequality and the fact that $\E[\exp(q^\t z)] \leq
\exp(\|q\|^2/2)$ for any vector $q \in \R^m$,
\begin{align*}
\E\left[ \exp\left( \alpha^\t \bigl(\sqrt{m} \Theta^\t e_i\bigr)
\right) \right]
& = \E\left[ \exp\left( \alpha^\t \bigl(\sqrt{m} v\bigr) \right) \right]
\\
& = \E\left[ \E\left[ \exp\left( \frac{\sqrt{m}}{\E[l]} \alpha^\t (\E[l]
v) \right) \ \Big| \ v \right] \right] \\
& \leq \E\left[ \exp\left( \frac{\sqrt{m}}{\E[l]} \alpha^\t (lv) \right)
\right] \\
& = \E\left[ \exp\left( \frac{\sqrt{m}}{\E[l]} \alpha^\t z \right)
\right] \\
& \leq \exp\left( \frac{\|\alpha\|^2 m}{2 \E[l]^2} \right) \\
& \leq \exp\left( \frac{\|\alpha\|^2}{2} \left( 1 - \frac1{4m} - \frac1{360m^3}
\right)^{-2} \right)
\end{align*}
where the last inequality is due to the following lower estimate for $\chi$
random variables:
\[ \E[l] \geq \sqrt{m} \left( 1 - \frac1{4m} - \frac1{360m^3} \right) . \]
Therefore, the condition~\eqref{eq:rotation-cond} is satisfied with $\kappa
= 1 + O(1/m)$.
\end{example}
\begin{example} \label{example:hadamard}
Let $m$ be a power of two, and let $\Theta := H\diag(s)/\sqrt{m}$, where $H
\in \{\pm1\}^{m \times m}$ is the $m \times m$ Hadamard matrix, and $s :=
(s_1,s_2,\dotsc,s_n) \in \{\pm1\}^m$ is a vector of $m$ Rademacher
variables (\emph{i.e.}, $s_1,s_2,\dotsc,s_m$ are i.i.d.~with
$\Pr[s_1=1]=\Pr[s_1=-1]=1/2$).
It is easy to check that $\Theta$ is an orthogonal matrix.
The random rotation $\Theta$ is a key component of the fast
Johnson-Lindenstrauss transform of~\cite{AilCha09}, also used
by~\cite{DMMS10}.
It is especially important for the present application because it can be
applied to vectors with $O(m \log m)$ operations, which is significantly
faster than the $\Omega(m^2)$ running time of na\"ive matrix-vector
multiplication.

For each $i=1,2,\dotsc,m$, the distribution of $\sqrt{m} \Theta^\t e_i$ is
the same as that of $s$, and therefore
\[
\E\left[ \exp\left( \alpha^\t \bigl(\sqrt{m} \Theta^\t e_i\bigr)
\right) \right]
= \E\left[ \exp\left( \alpha^\t s \right) \right]
\leq \exp(\|\alpha\|^2 / 2)
\]
where the last step follows by Hoeffding's inequality.
Therefore, the condition~\eqref{eq:rotation-cond} is satisfied with $\kappa
= 1$.
\end{example}

\section{Proofs of Theorem~\ref{theorem:ols} and Theorem~\ref{theorem:ridge}}
\label{section:ridge-proof}

The proof of Theorem~\ref{theorem:ridge} uses the decomposition of $\|\blh
- \beta\|_\Sig^2$ in Proposition~\ref{proposition:decomposition}, and then
bounds each term using the lemmas proved in this section.

The proof of Theorem~\ref{theorem:ols} omits one term from the
decomposition in Proposition~\ref{proposition:decomposition} due to the
fact that $\beta = \bl$ when $\lambda = 0$; and it uses a slightly simpler
argument to handle the effect of noise (Lemma~\ref{lemma:noise-ols} rather
than Lemma~\ref{lemma:noise}), which reduces the number of lower-order
terms.
Other than these differences, the proof is the same as that for
Theorem~\ref{theorem:ridge} in the special case of $\lambda = 0$.

Define
\begin{align}
\Sigl & := \Sig + \lambda I , \label{eq:Sigl} \\
\Sighl & := \Sigh + \lambda I , \label{eq:Sighl} \quad \text{and} \\
\Deltal & := \Sigl^{-1/2} (\Sigh - \Sig) \Sigl^{-1/2} \label{eq:Deltal} \\
& \hphantom{:}= \Sigl^{-1/2} (\Sighl - \Sigl) \Sigl^{-1/2} . \nonumber
\end{align}

Recall the basic decomposition from
Proposition~\ref{proposition:decomposition}:
\[
\|\blh - \beta\|_\Sig^2
\leq
\left( \|\bl - \beta\|_\Sig
+ \|\blb - \bl\|_\Sig
+ \|\blh - \blb\|_\Sig \right)^2
.
\]
Section~\ref{section:basic-properties} first establishes basic properties
of $\beta$ and $\bl$, which are then used to bound $\|\bl -
\beta\|_\Sig^2$; this part is exactly the same as the standard fixed design
analysis of ridge regression.
Section~\ref{section:matrix-error} employs probability tail inequalities
for the spectral and Frobenius norms of random matrices to bound the matrix
errors in estimating $\Sig$ with $\Sigh$.
Finally, Section~\ref{section:approximation-error} and
Section~\ref{section:noise} bound the contributions of approximation error
(in $\|\blb - \bl\|_\Sig^2$) and noise (in $\|\blh - \blb\|_\Sig^2$),
respectively, using probability tail inequalities for random vectors as
well as the matrix error bounds for $\Sigh$.

\subsection{Basic properties of $\beta$ and $\bl$, and the effect of
regularization}
\label{section:basic-properties}

The following propositions are well known in the study of inverse problems:

\begin{proposition}[Normal equations] \label{proposition:normal}
$\E[\dotp{w,x} y] = \E[\dotp{w,x} \dotp{\beta,x}]$ for any $w$.
\end{proposition}
\begin{proof}
It suffices to prove the claim for $w = v_j$.
Since $\E[\dotp{v_j,x} \dotp{v_{j'},x}] = 0$ for $j' \neq j$, it follows
that $\E[\dotp{v_j,x} \dotp{\beta,x}]
= \sum_{j'} \beta_{j'} \E[\dotp{v_j,x} \dotp{v_{j'},x}]
= \beta_j \E[\dotp{v_j,x}^2]
= \E[\dotp{v_j,x} y]$, where the last equality follows from the definition
of $\beta$ in~\eqref{eq:beta}.
\end{proof}

\begin{proposition}[Excess mean squared error] \label{proposition:regret}
$\E[(\dotp{w,x}-y)^2] - \E[(\dotp{\beta,x}-y)^2]
= \E[\dotp{w-\beta,x}^2]$ for any $w$.
\end{proposition}
\begin{proof}
Directly expanding the squares in the expectations reveals that
\begin{align*}
\lefteqn{
\E[(\dotp{w,x}-y)^2] - \E[(\dotp{\beta,x}-y)^2]
} \\
& = \E[\dotp{w,x}^2]
- 2\E[\dotp{w,x}y]
+ 2\E[\dotp{\beta,x}y]
- \E[\dotp{\beta,x}^2]
\\
& = \E[\dotp{w,x}^2]
- 2\E[\dotp{w,x}\dotp{\beta,x}]
+ 2\E[\dotp{\beta,x}\dotp{\beta,x}]
- \E[\dotp{\beta,x}^2]
\\
& = \E[\dotp{w,x}^2 - 2\dotp{w,x}\dotp{\beta,x} + \dotp{\beta,x}^2]
\\
& = \E[\dotp{w-\beta,x}^2]
\end{align*}
where the third equality follows from Proposition~\ref{proposition:normal}.
\end{proof}

\begin{proposition}[Shrinkage] \label{proposition:shrinkage}
For any $j$,
\[ \dotp{v_j,\bl} = \frac{\lambda_j}{\lambda_j + \lambda} \beta_j . \]
\end{proposition}
\begin{proof}
Since $(\Sig + \lambda I)^{-1} = \sum_j (\lambda_j + \lambda)^{-1} v_j
\otimes v_j$,
\[
\dotp{v_j,\bl}
= \dotp{v_j,(\Sig + \lambda I)^{-1} \E[xy]}
= \frac{1}{\lambda_j + \lambda} \E[\dotp{v_j,x}y]
= \frac{\lambda_j}{\lambda_j + \lambda}
\frac{\E[\dotp{v_j,x}y]}{\dotp{v_j,x}^2}
= \frac{\lambda_j}{\lambda_j + \lambda} \beta_j
.
\]
\end{proof}

\begin{proposition}[Effect of regularization] \label{proposition:regularization}
\[ \|\beta-\bl\|_\Sig^2
= \sum_j \frac{\lambda_j}{(\frac{\lambda_j}{\lambda} + 1)^2} \beta_j^2 .
\]
\end{proposition}
\begin{proof}
By Proposition~\ref{proposition:shrinkage},
\[ \dotp{v_j,\beta-\bl}
= \beta_j - \frac{\lambda_j}{\lambda_j + \lambda} \beta_j
= \frac{\lambda}{\lambda_j + \lambda} \beta_j
. \]
Therefore,
\[
\|\beta-\bl\|_\Sig^2
= \sum_j \lambda_j
\left( \frac{\lambda}{\lambda_j + \lambda} \beta_j \right)^2
= \sum_j \frac{\lambda_j}{(\frac{\lambda_j}{\lambda} + 1)^2} \beta_j^2
.
\]
\end{proof}

\subsection{Effect of errors in $\Sigh$}
\label{section:matrix-error}

\begin{lemma}[Spectral norm error in $\Sigh$] \label{lemma:deltal-spectral}
Assume Condition~\ref{cond:leverage}
(with parameter $\rhol$) holds.
Pick $t > \max\{0, 2.6 - \log\dlt\}$.
With probability at least $1-e^{-t}$,
\[
\|\Deltal\|
\leq \sqrt{\frac{4\rhol^2\dl(\log\dlt + t)}{n}} + \frac{2\rhol^2\dl
(\log\dlt + t)}{3n}
\]
where $\Deltal$ is defined in~\eqref{eq:Deltal}.
\end{lemma}
\begin{proof}
The claim is a consequence of the tail inequality from
Lemma~\ref{lemma:matrix-bernstein}.
First, define
\[ \xt := \Sigl^{-1/2} x
\quad\text{and}\quad
\Sigt := \Sigl^{-1/2} \Sig \Sigl^{-1/2}
\]
(where $\Sigl$ is defined in~\eqref{eq:Sigl}),
and let
\begin{align*}
Z & := \xt \otimes \xt - \Sigt \\
& \hphantom{:}= \Sigl^{-1/2} (x \otimes x - \Sig) \Sigl^{-1/2}
\end{align*}
so $\Deltal = \Eh[Z]$.
Observe that $\E[Z] = 0$ and
\[
\|Z\|
= \max\{ \lambda_{\max}[Z], \lambda_{\max}[-Z] \}
\leq \max\{ \|\xt\|^2, 1 \}
\leq \rhol^2 \dl
\]
where the second inequality follows from Condition~\ref{cond:leverage}.
Moreover,
\[
\E[Z^2]
= \E[(\xt \otimes \xt)^2] - \Sigt^2
= \E[\|\xt\|^2 (\xt \otimes \xt)] - \Sigt^2
\]
so
\begin{align*}
\lambda_{\max}[\E[Z^2]]
& \leq \lambda_{\max}[\E[(\xt \otimes \xt)^2]]
\leq \rhol^2 \dl \lambda_{\max}[\Sigt]
\leq \rhol^2 \dl \\
\tr(\E[Z^2])
& \leq \tr(\E[\|\xt\|^2 (\xt \otimes \xt)])
\leq \rhol^2 \dl \tr(\Sigt)
= \rhol^2 \dl^2
.
\end{align*}
The claim now follows from Lemma~\ref{lemma:matrix-bernstein} (recall that
$\dlt = \max\{ 1, \dl \}$).
\end{proof}

\begin{lemma}[Relative spectral norm error in $\Sighl$]
\label{lemma:deltal-ratio}
If $\|\Deltal\| < 1$ where $\Deltal$ is defined in~\eqref{eq:Deltal},
then
\[ \|\Sigl^{1/2} \Sighl^{-1} \Sigl^{1/2}\| \leq \frac{1}{1 - \|\Deltal\|}
\]
where $\Sigl$ is defined in~\eqref{eq:Sigl} and $\Sighl$ is defined
in~\eqref{eq:Sighl}.
\end{lemma}
\begin{proof}
Observe that
\begin{align*}
\Sigl^{-1/2} \Sighl \Sigl^{-1/2}
& = \Sigl^{-1/2} (\Sigl + \Sighl - \Sigl) \Sigl^{-1/2} \\
& = I + \Sigl^{-1/2} (\Sighl - \Sigl) \Sigl^{-1/2} \\
& = I + \Deltal
,
\end{align*}
and that
\[ \lambda_{\min}[I + \Deltal] \geq 1 - \|\Deltal\| > 0 \]
by the assumption $\|\Deltal\| < 1$ and Weyl's theorem~\cite{Horn:1985:MA}.
Therefore,
\[
\|\Sigl^{1/2} \Sighl^{-1} \Sigl^{1/2}\|
= \lambda_{\max}[(\Sigl^{-1/2} \Sighl \Sigl^{-1/2})^{-1}]
= \lambda_{\max}[(I + \Deltal)^{-1}]
= \frac{1}{\lambda_{\min}[I + \Deltal]}
\leq \frac{1}{1 - \|\Delta\|}
.
\]
\end{proof}

\begin{lemma}[Frobenius norm error in $\Sigh$] \label{lemma:deltal-frobenius}
Assume Condition~\ref{cond:leverage}
(with parameter $\rhol$) holds.
Pick any $t > 0$.
With probability at least $1-e^{-t}$,
\begin{align*}
\|\Deltal\|_\F
& \leq \sqrt{\frac{\E[\|\Sigl^{-1/2}x\|^4] - \dll}{n}} (1 + \sqrt{8t})
+ \frac{4\sqrt{\rhol^4\dl^2 + \dll}t}{3n}
\\
& \leq \sqrt{\frac{\rhol^2\dl^2 - \dll}{n}} (1 + \sqrt{8t})
+ \frac{4\sqrt{\rhol^4\dl^2 + \dll}t}{3n}
\end{align*}
where $\Deltal$ is defined in~\eqref{eq:Deltal}.
\end{lemma}
\begin{proof}
The claim is a consequence of the tail inequality in
Lemma~\ref{lemma:vector-bernstein}.
As in the proof of Lemma~\ref{lemma:deltal-spectral}, define $\xt :=
\Sigl^{-1/2} x$ and $\Sigt := \Sigl^{-1/2} \Sig \Sigl^{-1/2}$, and let $Z
:= \xt \otimes \xt - \Sigt$ so $\Deltal = \Eh[Z]$.
Now endow the space of self-adjoint linear operators with the inner product
given by $\dotp{A,B}_\F := \tr(AB)$, and note that this inner product
induces the Frobenius norm $\|M\|_\F = \dotp{M,M}_\F$.
Observe that $\E[Z] = 0$ and
\begin{align*}
\|Z\|_\F^2
& = \dotp{ \xt \otimes \xt - \Sigt, \xt \otimes \xt - \Sigt }_\F
\\
& =
\dotp{ \xt \otimes \xt, \xt \otimes \xt }_\F
- 2\dotp{ \xt \otimes \xt, \Sigt }_\F
+ \dotp{ \Sigt, \Sigt }_\F
\\
& = \|\xt\|^4 - 2\|\xt\|_{\Sigt}^2 + \tr(\Sigt^2)
\\
& = \|\xt\|^4 - 2\|\xt\|_{\Sigt}^2 + \dll
\\
& \leq \rhol^4 \dl^2 + \dll ,
\end{align*}
where the inequality follows from Condition~\ref{cond:leverage}.
Moreover,
\begin{align*}
\E[\|Z\|_\F^2]
& = 
\E[\dotp{ \xt \otimes \xt, \xt \otimes \xt }_\F]
- \dotp{ \Sigt, \Sigt }_\F
\\
& = \E[\|\xt\|^4] - \dll
\\
& \leq \rhol^2 \dl \E[\|\xt\|^2] - \dll
\\
& = \rhol^2 \dl^2 - \dll ,
\end{align*}
where the inequality again uses Condition~\ref{cond:leverage}.
The claim now follows from Lemma~\ref{lemma:vector-bernstein}.
\end{proof}

\subsection{Effect of approximation error}
\label{section:approximation-error}

\begin{lemma}[Effect of approximation error] \label{lemma:approximation}
Assume Condition~\ref{cond:leverage} (with
parameter $\rhol$) and Condition~\ref{cond:bias} (with parameter
$\bbiasl$) hold.
Pick any $t > 0$.
If $\|\Deltal\| < 1$ where $\Deltal$ is defined in~\eqref{eq:Deltal},
then
\[
\|\blb - \bl\|_\Sig
\leq
\frac{1}{1 - \|\Deltal\|}
\| \Eh[x \biasl(x) - \lambda \bl] \|_{\Sigl^{-1}} ,
\]
where
$\blb$ is defined in~\eqref{eq:blb},
$\bl$ is defined in~\eqref{eq:bl},
$\biasl(x)$ is defined in~\eqref{eq:biasl},
and $\Sigl$ is defined in~\eqref{eq:Sigl}.
Moreover, with probability at least $1-e^{-t}$,
\begin{align*}
\lefteqn{
\| \Eh[x \biasl(x) - \lambda \bl] \|_{\Sigl^{-1}}
} \\
& \leq
\sqrt{
\frac{\E[\| \Sigl^{-1/2} (x \biasl(x) - \lambda \bl)\|^2]}{n}}
(1 + \sqrt{8t})
+ \frac{4(\bbiasl \sqrt{\dl} + \|\beta-\bl\|_\Sig)t}{3n}
\\
& \leq
\sqrt{\frac{2(\rhol^2 \dl \E[\biasl(x)^2] + \|\beta-\bl\|_\Sig^2)}{n}}
(1 + \sqrt{8t})
+ \frac{4(\bbiasl \sqrt{\dl} + \|\beta-\bl\|_\Sig)t}{3n}
.
\end{align*}
\end{lemma}
\begin{proof}
By the definitions of $\blb$ and $\bl$,
\begin{align*}
\blb - \bl
& = \Sighl^{-1} \left( \Eh[x\E[y|x]] - \Sighl \bl \right) \\
& =
\Sigl^{-1/2} (\Sigl^{1/2} \Sighl^{-1} \Sigl^{1/2}) \Sigl^{-1/2}
\left( \Eh[x(\bias(x) + \dotp{\beta,x})] - \Sigh \bl - \lambda \bl \right)
\\
& =
\Sigl^{-1/2} (\Sigl^{1/2} \Sighl^{-1} \Sigl^{1/2}) \Sigl^{-1/2}
\left( \Eh[x(\bias(x) + \dotp{\beta,x} - \dotp{\bl,x})] - \lambda \bl
\right)
\\
& =
\Sigl^{-1/2} (\Sigl^{1/2} \Sighl^{-1} \Sigl^{1/2}) \Sigl^{-1/2}
\left( \Eh[x \biasl(x) - \lambda \bl] \right)
.
\end{align*}
Therefore, using the submultiplicative property of the spectral norm,
\begin{align*}
\|\blb - \bl\|_\Sig
& \leq
\|\Sig^{1/2} \Sigl^{-1/2}\|
\|\Sigl^{1/2} \Sighl^{-1} \Sigl^{1/2}\|
\| \Eh[x \biasl(x) - \lambda \bl] \|_{\Sigl^{-1}}
\\
& \leq
\frac{1}{1 - \|\Deltal\|}
\| \Eh[x \biasl(x) - \lambda \bl] \|_{\Sigl^{-1}}
\end{align*}
where the second inequality follows from Lemma~\ref{lemma:deltal-ratio} and
because
\[
\|\Sig^{1/2} \Sigl^{-1/2}\|^2
= \lambda_{\max}[\Sigl^{-1/2} \Sig \Sigl^{-1/2}]
= \max_i \frac{\lambda_i}{\lambda_i + \lambda}
\leq 1
.
\]

The second part of the claim is a consequence of the tail inequality in
Lemma~\ref{lemma:vector-bernstein}.
Observe that $\E[x\bias(x)] = \E[x(\E[y|x] - \dotp{\beta,x})] = 0$ by
Proposition~\ref{proposition:normal}, and that $\E[x\dotp{\beta - \bl,x}] -
\lambda \bl = \Sig\beta - (\Sig + \lambda I) \bl = 0$.
Therefore,
\[ \E[\Sigl^{-1/2}(x\biasl(x) - \lambda\bl)]
= \Sigl^{-1/2} \E[x(\bias(x) + \dotp{\beta - \bl,x}) - \lambda\bl]
= 0
. \]
Moreover, by Proposition~\ref{proposition:shrinkage} and
Proposition~\ref{proposition:regularization},
\begin{align}
\|\lambda \Sigl^{-1/2} \bl\|^2
& = \sum_j \frac{\lambda^2}{\lambda_j + \lambda} \dotp{v_j,\bl}^2
\nonumber \\
& = \sum_j \frac{\lambda^2}{\lambda_j + \lambda} \left(
\frac{\lambda_j}{\lambda_j + \lambda} \beta_j \right)^2
\nonumber \\
& \leq \sum_j \frac{\lambda^2}{\lambda_j + \lambda} \left(
\frac{\lambda_j}{\lambda_j + \lambda} \right) \beta_j^2
\nonumber \\
& = \sum_j \frac{\lambda_j}{(\frac{\lambda_j}{\lambda} + 1)^2} \beta_j^2
\nonumber \\
& = \|\beta - \bl\|_\Sig^2
\label{eq:biasl-lb}
.
\end{align}
Combining the inequality from~\eqref{eq:biasl-lb} with
Condition~\ref{cond:bias} and the triangle inequality, it follows that
\begin{align*}
\|\Sigl^{-1/2} (x\biasl(x) - \lambda \bl)\|
& \leq
\|\Sigl^{-1/2} x\biasl(x)\|
+ \|\lambda \Sigl^{-1/2} \bl\| \\
& \leq \bbiasl \sqrt{\dl} + \|\beta - \bl\|_\Sig
.
\end{align*}
Finally, by the triangle inequality, the fact $(a+b)^2 \leq 2(a^2+b^2)$,
the inequality from~\eqref{eq:biasl-lb}, and Condition~\ref{cond:leverage},
\begin{align*}
\E[\|\Sigl^{-1/2} (x\biasl(x) - \lambda \bl)\|^2]
& \leq 2 (\E[\|\Sigl^{-1/2}x\biasl(x)\|^2] + \|\bl-\beta\|_\Sig^2) \\
& \leq 2 (\rhol^2 \dl \E[\biasl(x)^2] + \|\bl-\beta\|_\Sig^2)
.
\end{align*}
The claim now follows from Lemma~\ref{lemma:vector-bernstein}.
\end{proof}

\subsection{Effect of noise}
\label{section:noise}

\begin{lemma}[Effect of noise, $\lambda = 0$] \label{lemma:noise-ols}
Assume 
$\lambda = 0$.
Assume Condition~\ref{cond:noise} (with parameter $\snoise$) holds.
Pick any $t > 0$.
With probability at least $1-e^{-t}$,
either $\|\Delta_0\| \geq 1$, or
\[
\|\Delta_0\| < 1 \quad\text{and}\quad
\|\bb_0 - \bh_0\|_\Sig^2 \leq
\frac1{1-\|\Delta_0\|} \cdot \frac{\snoise^2 ( d + 2\sqrt{dt} + 2 t)}{n}
,
\]
where $\Delta_0$ is defined in~\eqref{eq:Deltal}.
\end{lemma}
\begin{proof}
Observe that
\[ \|\bb_0 - \bh_0\|_\Sig^2
\leq \|\Sig^{1/2} \Sigh^{-1/2}\|^2
\|\bb_0 - \bh_0\|_{\Sigh}^2
= \|\Sig^{1/2} \Sigh^{-1} \Sig^{1/2}\|
\|\bb_0 - \bh_0\|_{\Sigh}^2
;
\]
and if $\|\Delta_0\| < 1$, then $\|\Sig^{1/2} \Sigh^{-1} \Sig^{1/2}\| \leq
1/(1 - \|\Delta_0\|)$ by Lemma~\ref{lemma:deltal-ratio}.

Let $\xi := (\noise(x_1),\noise(x_2),\dotsc,\noise(x_n))$ be the random
vector whose $i$-th component is $\noise(x_i) = y_i - \E[y_i|x_i]$.
By the definition of $\bh_0$ and $\bb_0$
\[
\|\bh_0 - \bb_0\|_{\Sigh}^2
= \|\Sigh^{-1/2} \Eh[x(y - \E[y|x])]\|^2
= \xi^\t \wh{K} \xi ,
\]
where $\wh{K} \in \R^{n \times n}$ is the symmetric matrix whose $(i,j)$-th
entry is $\wh{K}_{i,j} := n^{-2} \dotp{\Sigh^{-1/2} x_i, \Sigh^{-1/2}
x_j}$.
Note that the nonzero eigenvalues of $\wh{K}$ are the same as those of
\[ \frac1n \Eh\left[ (\Sigh^{-1/2} x) \otimes (\Sigh^{-1/2} x) \right]
= \frac1n \Sigh^{-1/2} \Sigh \Sigh^{-1/2}
= \frac1n I
.
\]
By Lemma~\ref{lemma:quadratic}, with probability at least $1-e^{-t}$
(conditioned on $x_1,x_2,\dotsc,x_n$),
\[
\xi^\t \wh{K} \xi \leq \snoise^2 ( \tr(\wh{K}) + 2\sqrt{\tr(\wh{K}^2)t} +
2\lambda_{\max}(\wh{K}) t) \\
= \frac{\snoise^2 (d + 2\sqrt{dt} + 2t)}{n}
.
\]
The claim follows.
\end{proof}

\begin{lemma}[Effect of noise, $\lambda \geq 0$] \label{lemma:noise}
Assume Condition~\ref{cond:noise} (with
parameter $\snoise$) holds.
Pick any $t > 0$.
Let $K$ be the $n \times n$ symmetric matrix whose $(i,j)$-th entry is
\[ K_{i,j} := \frac{1}{n^2} \dotp{\Sig^{1/2} \Sighl^{-1} x_i,
\Sig^{1/2} \Sighl^{-1} x_j} ,
\]
where $\Sighl$ is defined in~\eqref{eq:Sighl}.
With probability at least $1-e^{-t}$,
\[
\|\blb - \blh\|_\Sig^2
\leq \snoise^2 ( \tr(K) + 2\sqrt{\tr(K)\lambda_{\max}(K)t} +
2\lambda_{\max}(K) t)
.
\]
Moreover, if $\|\Deltal\| < 1$ where $\Deltal$ is defined
in~\eqref{eq:Deltal},
then
\[
\lambda_{\max}(K) \leq \frac{1}{n(1-\|\Deltal\|)}
\quad\text{and}\quad
\tr(K) \leq \frac{\dll + \sqrt{\dll \|\Deltal\|_\F^2}}{n(1-\|\Deltal\|)^2}
.
\]
\end{lemma}
\begin{proof}
Let $\xi := (\noise(x_1),\noise(x_2),\dotsc,\noise(x_n))$ be the random
vector whose $i$-th component is $\noise(x_i) = y_i - \E[y_i|x_i]$.
By the definition of $\blh$, $\blb$, and $K$,
\[
\|\blh - \blb\|_\Sig^2
= \|\Sighl^{-1} \Eh[x(y - \E[y|x])]\|_\Sig^2
= \xi^\t K \xi
.
\]
By Lemma~\ref{lemma:quadratic}, with probability at least $1-e^{-t}$
(conditioned on $x_1,x_2,\dotsc,x_n$),
\begin{align*}
\xi^\t K \xi
& \leq \snoise^2 ( \tr(K) + 2\sqrt{\tr(K^2)t} + 2\lambda_{\max}(K) t) \\
& \leq \snoise^2 ( \tr(K) + 2\sqrt{\tr(K)\lambda_{\max}(K)t} +
2\lambda_{\max}(K) t) ,
\end{align*}
where the second inequality follows from von Neumann's theorem~\cite{Horn:1985:MA}.

Note that the nonzero eigenvalues of $K$ are the same as that of
\[ \frac{1}{n} \Eh\left[ (\Sig^{1/2} \Sighl^{-1} x) \otimes (\Sig^{1/2}
\Sighl^{-1} x) \right]
= \frac{1}{n} \Sig^{1/2} \Sighl^{-1} \Sigh \Sighl^{-1} \Sig^{1/2}
. \]
To bound $\lambda_{\max}(K)$, observe that by the submultiplicative
property of the spectral norm and Lemma~\ref{lemma:deltal-ratio},
\begin{align*}
n \lambda_{\max}(K)
& = \|\Sig^{1/2} \Sighl^{-1} \Sigh^{1/2} \|^2 \\
& \leq
\|\Sig^{1/2} \Sigl^{-1/2}\|^2
\|\Sigl^{1/2} \Sighl^{-1/2}\|^2
\|\Sighl^{-1/2} \Sigh^{1/2}\|^2 \\
& \leq
\|\Sigl^{1/2} \Sighl^{-1/2}\|^2
\\
& =
\|\Sigl^{1/2} \Sighl^{-1} \Sigl^{1/2}\|
\\
& \leq \frac{1}{1 - \|\Deltal\|}
.
\end{align*}
To bound $\tr(K)$, first define the $\lambda$-whitened versions of $\Sig$,
$\Sigh$, and $\Sighl$ as
\begin{align*}
\Sigw & := \Sigl^{-1/2} \Sig \Sigl^{-1/2} , \\
\Sighw & := \Sigl^{-1/2} \Sigh \Sigl^{-1/2} , \\
\Sighlw & := \Sigl^{-1/2} \Sighl \Sigl^{-1/2}
.
\end{align*}
Using these definitions with the cycle property of the trace,
\begin{align*}
n \tr(K)
& = \tr(\Sig^{1/2} \Sighl^{-1} \Sigh \Sighl^{-1} \Sig^{1/2})
\\
& = \tr(\Sighl^{-1} \Sigh \Sighl^{-1} \Sig)
\\
& = \tr(\Sighlw^{-1} \Sighw \Sighlw^{-1} \Sigw)
.
\end{align*}
Let $\{ \lambda_j[M] \}$ denote the eigenvalues of a linear operator $M$.
By von Neumann's theorem~\cite{Horn:1985:MA},
\[ \tr(\Sighlw^{-1} \Sighw \Sighlw^{-1} \Sigw)
\leq \sum_j \lambda_j[\Sighlw^{-1} \Sighw \Sighlw^{-1}] \lambda_j[\Sigw]
\]
and by Ostrowski's theorem~\cite{Horn:1985:MA},
\[ \lambda_j[\Sighlw^{-1} \Sighw \Sighlw^{-1}]
\leq \lambda_{\max}[\Sighlw^{-2}] \lambda_j[\Sighw] . \]
Therefore
\begin{align*}
\tr(\Sighlw^{-1} \Sighw \Sighlw^{-1} \Sigw)
& \leq \lambda_{\max}[\Sighlw^{-2}]
\sum_j \lambda_j[\Sighw] \lambda_j[\Sigw]
\\
& \leq \frac{1}{(1 - \|\Deltal\|)^2}
\sum_j \lambda_j[\Sighw] \lambda_j[\Sigw]
\\
& = \frac{1}{(1 - \|\Deltal\|)^2}
\sum_j
\left(
\lambda_j[\Sigw]^2 + (\lambda_j[\Sighw] - \lambda_j[\Sigw]) \lambda_j[\Sigw]
\right)
\\
& \leq \frac{1}{(1 - \|\Deltal\|)^2}
\left(
\sum_j
\lambda_j[\Sigw]^2
+ \sqrt{\sum_j (\lambda_j[\Sighw] - \lambda_j[\Sigw])^2}
\sqrt{\sum_j \lambda_j[\Sigw]^2}
\right)
\\
& = \frac{1}{(1 - \|\Deltal\|)^2}
\left(
\dll
+ \sqrt{\sum_j (\lambda_j[\Sighw] - \lambda_j[\Sigw])^2}
\sqrt{\dll}
\right)
\\
& \leq \frac{1}{(1 - \|\Deltal\|)^2}
\left(
\dll
+ \|\Sighw - \Sigw\|_\F \sqrt{\dll}
\right)
\\
& = \frac{1}{(1 - \|\Deltal\|)^2}
\left( \dll + \|\Deltal\|_\F \sqrt{\dll} \right) ,
\end{align*}
where the second inequality follows from Lemma~\ref{lemma:deltal-ratio},
the third inequality follows from Cauchy-Schwarz, and the fourth inequality
follows from Mirsky's theorem~\cite{Stewart90}.
\end{proof}

\subsubsection*{Acknowledgements}
The authors thank Dean Foster, David McAllester, and Robert Stine for many insightful discussions.

\bibliography{ridge}
\bibliographystyle{plain}

\appendix

\section{Probability tail inequalities}
\label{appendix:tail}

The following probability tail inequalities are used in our analysis.
These specific inequalities were chosen in order to satisfy the general
conditions set up in Section~\ref{section:datamodel}; however, our analysis
can specialize or generalize with the availability of other tail
inequalities of these sorts.

The first tail inequality is for positive semidefinite quadratic forms of a
subgaussian random vector.
It generalizes a standard tail inequality for Gaussian random vectors based
on linear combinations of $\chi^2$ random variables~\cite{LauMas00}.

\begin{lemma}[Quadratic forms of a subgaussian random vector;
\cite{HKZ_vector}]
\label{lemma:quadratic}
Let $\xi$ be a random vector taking values in $\R^n$ such that for some
$c \geq 0$,
\[
\E[\exp(\dotp{u,\xi})]
\leq \exp(c \|u\|^2 / 2),
\quad \forall u \in \R^n
.
\]
For all symmetric positive semidefinite matrices $K \succeq 0$, and all
$t > 0$,
\[
\Pr\biggl[ \xi^\t K \xi
> c \Bigl( \tr(K) + 2\sqrt{\tr(K^2)t} +
2\|K\|t \Bigr) \biggr] \leq e^{-t}
.
\]
\end{lemma}

The next lemma is a tail inequality for sums of bounded random vectors; it
is a standard application of Bernstein's inequality.
\begin{lemma}[Vector Bernstein bound; see, \emph{e.g.},
\cite{HKZ_vector}]
\label{lemma:vector-bernstein}
Let $x_1,x_2,\dotsc,x_n$ be independent random vectors such that
\[
\sum_{i=1}^n \E[ \|x_i\|^2 ] \leq v
\quad \text{and} \quad
\|x_i\| \leq r
\]
for all $i=1,2,\dotsc,n$, almost surely.
Let $s := x_1 + x_2 + \dotsb + x_n$.
For all $t > 0$,
\[
\Pr\left[ \|s\| > \sqrt{v} (1 + \sqrt{8t}) + (4/3) r t
\right] \leq e^{-t}
\]
\end{lemma}

The last tail inequality concerns the spectral accuracy of an empirical
second moment matrix.
\begin{lemma}[Matrix Bernstein bound; \cite{HKZ_matrix}]
\label{lemma:matrix-bernstein}
Let $X$ be a random matrix, and $r > 0$, $v > 0$, and $k >
0$ be such that, almost surely,
\begin{gather*}
\E[X] = 0 , \quad
\lambda_{\max}[X] \leq r , \quad
\lambda_{\max}[\E[X^2]] \leq v , \quad
\tr(\E[X^2]) \leq v k
.
\end{gather*}
If $X_1,X_2,\dotsc,X_n$ are independent copies of $X$, then for any $t >
0$,
\begin{equation*}
\Pr\left[
\lambda_{\max}\left[ \frac1n \sum_{i=1}^n X_i \right]
> \sqrt{\frac{2vt}{n}} + \frac{rt}{3n}
\right]
\leq k t (e^t - t - 1)^{-1}
.
\end{equation*}
If $t \geq 2.6$, then $t (e^t - t - 1)^{-1} \leq e^{-t/2}$.
\end{lemma}

\end{document}